\documentclass[final]{siamltex} 

\usepackage{graphicx,color}
\usepackage{amssymb,amsmath}

\usepackage[colorlinks,bookmarksopen,bookmarksnumbered,citecolor=red,urlcolor=red]{hyperref}

\newtheorem{Thm}{Theorem}
\newtheorem{Lem}[Thm]{Lemma}
\newtheorem{rem}[Thm]{Remark}

\def\bm{\boldsymbol}
\def\mb{\mathbf}
\def\bxi{\bm{\xi}}
\def\veps{\varepsilon}
\def\half{\frac{1}{\,2\,}}
\def\wt{\widetilde}
\let\dsp=\displaystyle
\def\cal{\mathcal}
\def\ie{{\it i.e.}}
\def\eg{{\it e.g.}}

\title{Coarse-graining molecular dynamics models using an extended Galerkin projection}
\author{Xiantao Li
\thanks{Department of Mathematics, the Pennsylvania State University, University Park,
Pennsylvania, 16802, U.S.A. ({\tt xli@math.psu.edu})
{The work of Li is supported by National Natural Science
              Foundation grant DMS1016582.}}}

\begin{document}

\maketitle

\begin{abstract}
We present a new framework for coarse-graining molecular dynamics models for crystalline solids.
The reduction method is based on a Galerkin projection to a subspace, whose dimension is much smaller than that of the full atomistic model. The subspace is expanded by adding more coarse-grain variables near the interface between lattice defects and the surrounding regions. This effectively
minimizes reflection of phonons at the interface. In this approach,
there is no need to pre-compute the memory function in the generalized Langevin equations, a typical
 model of interface conditions. Moreover, the variational formulation preserves the stability 
of mechanical equilibria.     
\end{abstract}



\section{Introduction}

Molecular dynamics (MD) models have been playing an increasingly important role in the study
of micro-mechanical systems. The main difficulty in the direct implementation of the MD model is the
huge number of atomic degrees of freedom. Simulating a realistic system over a time scale of practical
interest usually requires an extremely expensive computation. Therefore, it is of great interest to
develop coarse-grained (CG) MD models in which only a small number of representative degrees of
freedom are involved. For systems close to equilibrium, a coarse-graining energy can be obtained by
integrating out the remaining degrees of freedom in the free energy. This is the idea behind
most direct CG models, such as the one by Rudd and Broughton \cite{RuBr05,RuBr98} for crystalline solids. The equations of motion is then expressed using the Hamilton's principle. Notice, however, these models are not derived directly from the MD model.
In particular, for systems far from equilibrium, \eg, propagation and interaction of shock waves, this procedure can not be applied.  There have also been attempts to derive CG models using homogenization methods \cite{ChFi06,ClCh06,FiCh07}. Such asymptotic approach offers a means to derive a continuum mechanics model, although such procedure can only be applied to perfect crystals. In contrast,  Seleson et al \cite{SePa09} considered the upscaling to peridynamics (PD) models \cite{Sillings00}.

Another approach that has gained great popularity is to directly combine MD with continuum mechanics models. The idea is to retain the atomic-level description near critical areas where the molecular trajectories are of interest, while in the surrounding region, the MD model is replaced by a continuum mechanics model,  often in the form of partial differential equations. At the interface, a coupling condition is imposed to combine the two models, see \cite{AbBrBe98,EHu01,EHu02,Rodney03,ChLe03a,WaKa04,XiBe04,QuShCuMi05,ToLi05,ELi05,LiE05,TaHoLi06,TaHoLi06b,OdPrRoBa06,WaLiShu08,citeulike:2898408,BaPaBoGuLe08,LiYaE09,KrMi10,FaKr11}. 
Despite the many proposed methods, quite a few important issues have not been fully addressed. Examples include 
systems at finite temperature,  dislocation dynamics problems where dislocations pass through the interface, and adaptive selection of the MD region. In addition, stability and error analysis for such coupling methods is open.

In this paper, we will consider the Galerkin projection procedure to derive CG models for crystalline
systems. The method is motivated by the Galerkin approximations of some Hamiltonian PDEs, especially
the second-order wave equations. Examples include the finite element approximations of the
linear second-order hyperbolic equation \cite{Baker76,BakDou76,Dup73,Geveci84} and nonlinear
elasto-dynamics equation \cite{Mak93}. In fact, as demonstrated in \cite{ArGr05,EMing:2007b,SePa09}, for perfect crystals the MD model at zero temperature
can be approximated by the elasto-dynamics models with the constitutive relation given by the Cauchy-Born rule. Within the conventional Galerkin framework, the coarse-grained MD model can be derived by projecting the equations of motion to the subspace spanned by the basis functions. Near defects, this projection leaves the MD equation unchanged, while in the surrounding region, the CG equations coincide with those in \cite{WaLi03}, obtained by a least square procedure. The CG model obtained from this conventional Galerkin projection is also similar to the model developed by Chen and her coworkers \cite{ChLe03a,ChLe03b,XiTuMcCh11}, which were derived from the Irving-Kirkwood 
formalism \cite{irving1950statistical}.  But in their model, the issue of phonon reflection has not been addressed.

As we will show, the conventional Galerkin method often leads to
significant phonon reflections at the interface. The reflection can be attributed to the discrepancy
between the dispersion relations in the atomistic and CG regions, as well as the varying mesh size. In fact,
when solving wave equations, non-uniform mesh often causes inter-grid reflections 
\cite{Baz78,CeBa83}.

To minimize the phonon reflection, Wagner, Liu and their coworkers \cite{WaLi03,TaHoLi06,TaHoLi06b} 
developed a bridging scale method in which the system is divided into a coarse scale region and a MD region.
In the coarse scale region, the displacement of the atoms consists of a coarse scale part, the evolution of which is governed by equations derived from the Galerkin projection, and a fine scale part, which can be eliminated by
linearizing the MD model and using Laplace transform. This leads to a coupling condition at 
the interface, expressed in the form of a generalized Langevin equation (GLE). The reduction procedure is closely related 
to the derivation of absorbing boundary conditions for MD models \cite{AdDo76,AdDo74,Tully80,WaKa04,KarpovYu2005,KaWa05,LiE06,Li08b}.  
In order to implement this coupling condition, the memory kernel function in the generalized Langevin equation has to be precomputed. This calculation is a highly nontrivial task, especially because the functions tend to have very slow decay and they depend on the geometry of the interface. This numerical issue has been addressed in \cite{EHu01,EHu02,LiYaE09}, where the memory function is approximated by
functions with compact support.   A similar approach to the BS method is the perfectly matched (PMM) multiscale method \cite{ToLi05,LiLiu06}, in which the MD equations are modified by adding damping terms and modifying the force constant matrices at the interface to prevent artificial reflection. This method is easier to implement. However, it is  not clear how this can be correctly formulated in the case of coupling methods. Ideally, one should only remove the waves modes that are not captured by the coarse scale variables.  But no such selective damping mechanism is present in the design of the PMM method.

A common open problem in the methods mentioned above is  the well-posedness of the coupled models.  By imposing additional conditions at the interface, the variational structure is lost. 
In fact, for an isolated atomistic system, the stability of the boundary conditions is already
nontrivial \cite{Li08a}.
We will show some examples in this paper that some truncation schemes may lead to unstable models.   In this paper, we propose a new approach to extend the conventional Galerkin method. 
The idea is to expand the approximation space while maintaining the variational formulation. We start with a conventional Galerkin method. 
This is followed by a few refinement steps. At each step, we keep the mesh structure intact. Rather, we define additional 
variables at each nodal point. Therefore there is no need to re-mesh frequently to enhance the overall accuracy. 
The coarse-grained MD model is derived by projecting the equations of motion to the extended subspace. 
Therefore the variational structure is kept. 
As the expansion continues, we are building a hierarchy of ODE systems,
which in the limit is equivalent to the full atomistic model provided that the initial set of 
basis functions are chosen properly.
 In particular, the conventional Galerkin approximation is at the top of the hierarchy.   It can be obtained by dropping all the additional variables. From this perspective,
the large error that usually arises at the interface can be attributed to the premature truncation of the hierarchal structure. This structure is reminiscent of the BBGKY approximations of the many-body Liouville equation in non-equilibrium statistical physics \cite{Balescu76,Cerc90}. Our expansion procedure is in the spirit of higher moment models, \eg, the Bennett model \cite{AgYuBa01}. By 
introducing such higher order models at the interface, we are able to minimize phonon reflections
and improve the accuracy of the CG model.

This Galerkin projection approach bears similarity to the Mori-Zwanzig (MZ) projection formalism, a well known reduction methodology  in nonlinear statistical physics \cite{Mori65,Mori1965b,Zwanzig73}. For systems out-of-equilibrium, The MZ formalism has been used to derive coarse-grained molecular dynamics models \cite{Munakata85,CuCe02,IzVo06,Li2009c}, and it has proven to be a very useful way of thinking about this type of problems. In the MZ formalism, one uses
a projection operator to separate out coarse-grain variables, and then derive the effective dynamics. 
In the works on Zwanzig \cite{Zwanzig73} and the recent work of Chorin and his coworkers \cite{ChHaKu02,ChHaKu00,ChKaKu99,ChKaKu98},  conditional expectation is used as the projection. 
This is convenient when the initial data are sampled from an equilibrium distribution, \eg, the canonical ensemble.  In principle, one can use this approach to derive a GLE as a CG model \cite{Li2009c}. But it suffers the same problem mentioned above -- the computation of the memory kernel function. In fact, the expression of the memory function involves the inversion of a large matrix, with dimension equal to the total number of degrees of freedom in the system. On the other 
hand, crude approximations, \eg, by Dirac delta functions, are too ad hoc, and often introduce large modeling error.

The rest of the paper is organized as follows. We first describe the conventional Galerkin
method, applied to the MD model. We then show a simple numerical test to demonstrate
the accuracy of this method. Section \ref{sec: ext} presents the main idea of extended
Galerkin method, followed by some tips for practical implementation. In section \ref{sec: num}
we give some results of the numerical experiments. 

\section{The Conventional Galerkin Projection Method}

\subsection{The mathematical formulation}

\begin{figure}[htbp]
\begin{center}
\includegraphics[scale=0.7]{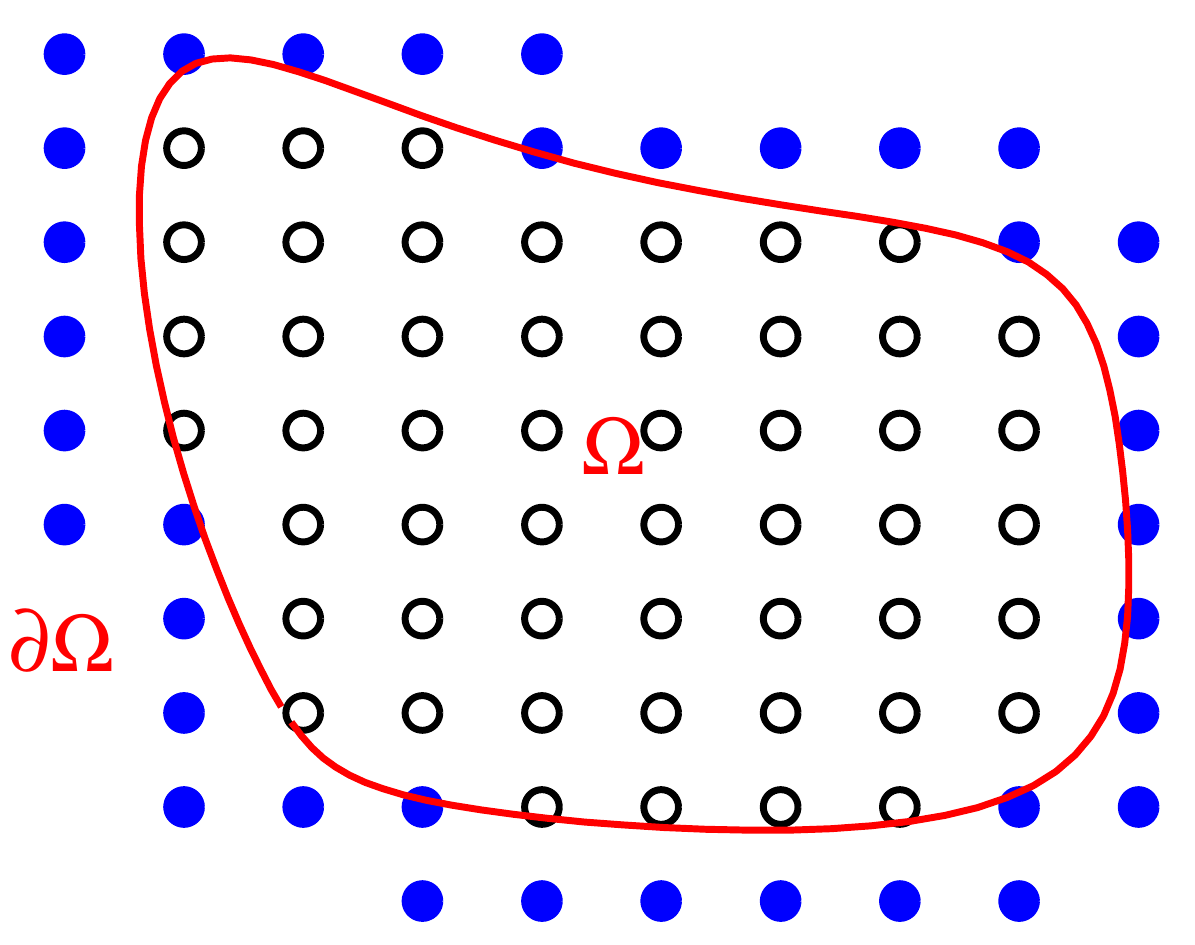}
\caption{A schematic illustration of the domain. Open circles: the atoms in the interior; Closed circles: the atoms at the boundary. }
\label{fig: domain}
\end{center}
\end{figure}

The formulation has strong analogy to the Galerkin projection of second order wave equations. But there 
are also crucial differences that will be emphasized. We consider a system of atoms with reference position in domain \(\Omega\) with boundary \(\partial \Omega\).  A schematic illustration is shown in Fig. \ref{fig: domain}.

We let \(\mb y_i\) and \(\mb x_i\) be respectively the reference and current position of the \(i\)th atom, with its displacement given by \(\mb u_i= \mb x_i - \mb y_i\).  For the molecular dynamics model, the atoms obey
the Newton's second law,
\begin{equation}\label{eq: md}
 m_i \ddot{\mb x}_i= -\nabla_{\mb x_i} V.
\end{equation}
Here \(m_i\) is the mass of the \(i\)th atom, and $\ddot{\mb x}_i$ denotes the second derivative in time.
Further \(V=V(\mb x_1, \mb x_2, \cdots, \mb x_N)\) is the interatomic potential, which will be assumed to be an empirical model. The boundary condition
is assumed to be of Dirichlet type, \ie, the position of the atoms in \(\partial \Omega\) is prescribed.

We first write the equations of motion in the following compact form,
\begin{equation}\label{eq: md0}
  \ddot{\mb u} = \mb f(\mb u),\;\; \mb f= -\nabla V.
\end{equation}
Here we have chosen to work with the displacement and grouped the displacement of all the atoms into one vector. We have set the mass to unity  and neglected possible external forces for simplicity 
as they will not significantly change the formulation.

To explain the numerical procedure, it is often helpful to separate out the harmonic and anharmonic parts of the interaction as follows,
\begin{equation}
  \mb f(\mb u)= -A \mb u + \mb f_N(\mb u).
\end{equation} 
Here \(\mb f_N\) indicates the anharmonic part of the force. The harmonic interaction can often be defined by a linearization around the reference configuration. See \cite{YaOrBh07} for such examples.  We will refer to the model
\begin{equation}\label{eq: md1}
  \ddot{\mb u} = -A \mb u,
\end{equation}
as the linearized MD model.

\begin{figure}[htbp]
\label{fig: interp}
\begin{center}
\includegraphics[scale=0.4]{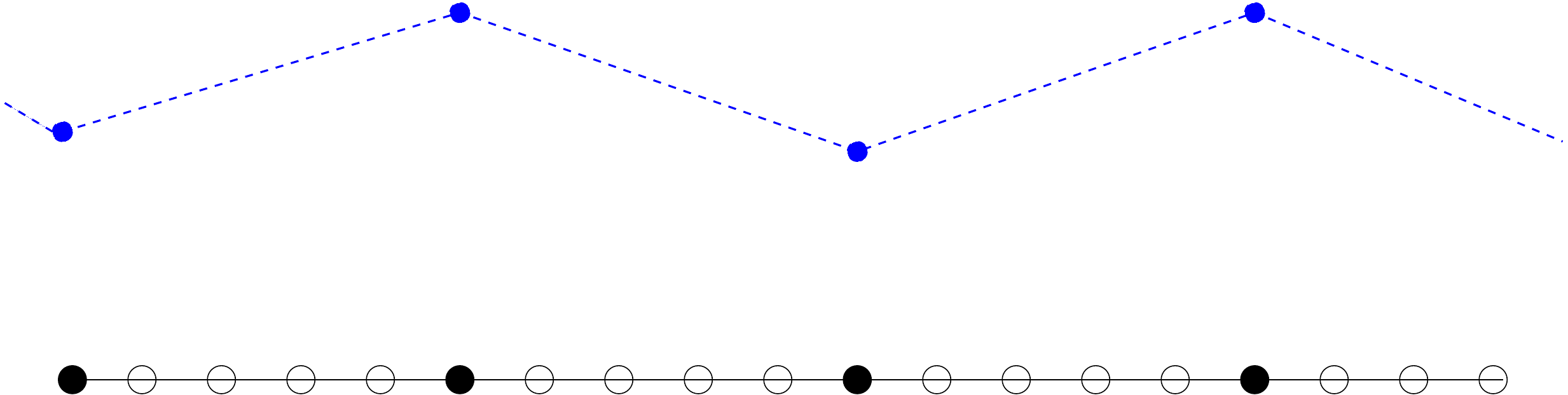}\\
\includegraphics[scale=0.4]{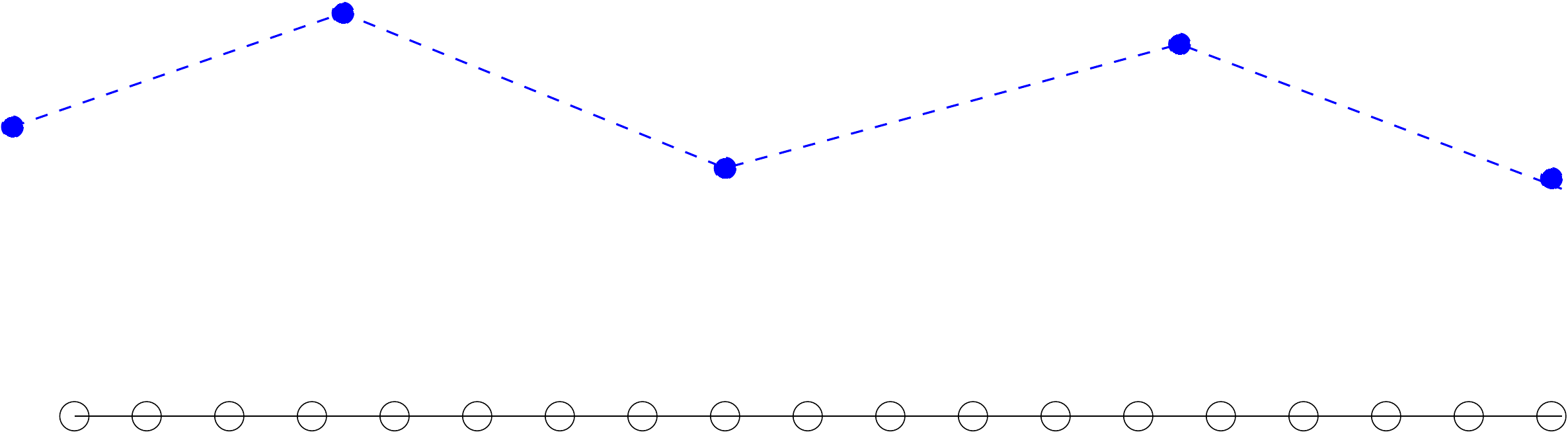}\\
\includegraphics[scale=0.4]{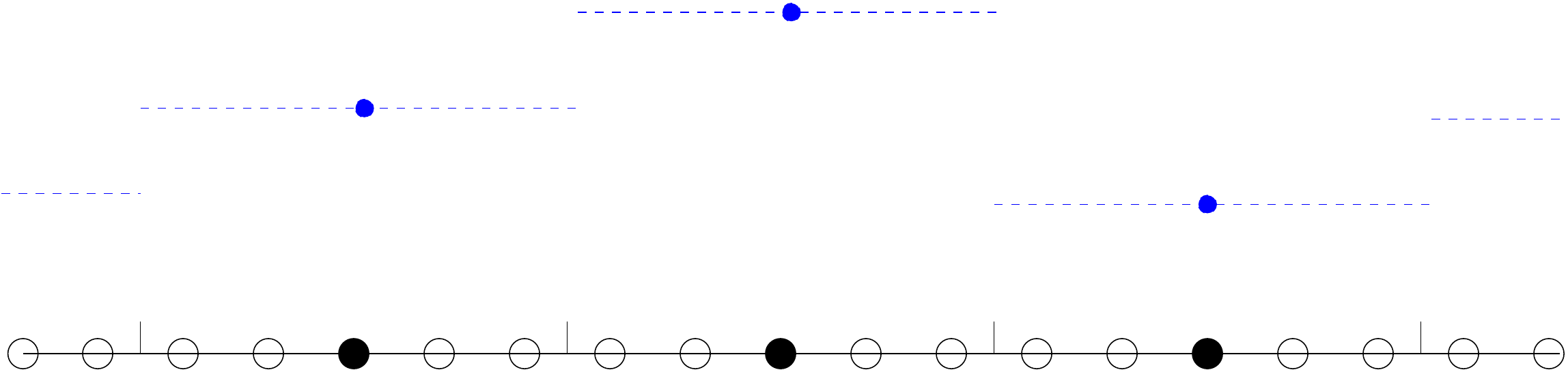}
\end{center}
\caption{Interpolation of a discrete function. Top: A piecewise linear interpolation in which the
nodal points coincide with the reference position of some of the atoms; Middle: Another piecewise linear interpolation 
where some of the nodal points are away from the atomic position; Bottom: A piecewise constant interpolation.}
\end{figure}

Next, let \(n\) be the dimension of \(\mb u\), and let \(X=\mathbb{R}^n\) be the entire space. To construct a space for
the approximate solutions, we choose a set of \(m\) independent nodal basis vectors, \(\big\{\varphi_i,1\le i\le m\big\}\). Let  \(\Phi=\Big[\varphi_1,\; \varphi_2,\; \cdots, \varphi_m  \Big]\) be a matrix, each column of which is a basis vector. We denote \(Y=\mathrm{Range}(\Phi)\) the \(m\)-dimensional subspace.  
The matrix \(\Phi\) may be constructed from an interpolation procedure for continuous problems, and each function \(\varphi_i: \Omega\to \mathbb{R}\) is a shape function.  One can also think of the matrix \(\Phi\) as 
an interpolation operator from \(\mathbb{R}^m\) to \(X\), and similarly, \(\Phi^T\) can be considered as a restriction operator, similar to the terminology  used in algebraic multi grid methods \cite{McCormick87}.   To better illustrate this idea, we give two examples. The first example is similar to a finite element representation with
piecewise linear functions. The top figure in Fig \ref{fig: interp} shows one such function in the subspace. In this case, the nodal points coincide with the reference position of some atoms. The middle figure shows a similar function, but the nodal points are away from atom positions. For the second example, we consider a piecewise constant function, which is similar to first order finite volume representation. We remark that the basis functions 
are {\it discrete} functions, although they can often be obtained by confining a continuous basis function to atomic positions. Near lattice defects, we often retain all the atoms in a neighborhood. For these atoms, we choose the trivial basis function which equals to one at an atomic position and zero at all other atoms. This part of the construction is due to the discreteness of the MD model and has no analogy in Galerkin approximations of PDEs.

\smallskip

 In \(X\), we define the usual inner product,
\begin{equation}
  \big(\mb u, \mb v \big)= \sum_{i} \mb u_i \cdot \mb  v_i.
\end{equation}
In the Hilbert space, we define the orthogonal projection to the subspace \(Y\), which in matrix notations, can be expressed as
\begin{equation}
  {P}= \Phi \big(\Phi^T \Phi)^{-1} \Phi^T.
\end{equation} In addition, we let \({Q}=I-{P}\) be the complementary projection. 

In addition to the inner product, we define,
\begin{equation}\label{eq: nonlinear}
 a(\mb u, \mb v)= -\big(\mb f(\mb u), \mb v).
\end{equation}
For linear models, this is reduced to,
\begin{equation}
a(\mb u, \mb v)= \mb u^T A  \mb v.
\end{equation}
Using some symmetry properties of \(A\), this quantity can often be written as \cite{AsMe76},
\[
a(\mb u, \mb v)= \sum_{i,j} (\mb u_i - \mb u_j)^T A_{ij}  (\mb v_i - \mb v_j).\]
Where \(A_{ij}\) is a force constant matrix.  This quadratic  form is similar to the bilinear form of elliptic PDEs.

We first discuss the conventional Galerkin approximation. Such an approach seeks an approximation \(\widehat{\mb u}: [0, T] \to Y\), such that,
 \(\ddot{\widehat{ \mb u}} \in L^2(0,T;X)\), and  for any test function \(\mb v \in Y\) that vanishes at the boundary,  there holds,
\begin{equation}\label{eq: Gal0}
  \big( \ddot{\widehat{\mb u}}, \mb v\big)  + a(\widehat{\mb u}, \mb v)=0.
\end{equation}

If we let \( {M}= \Phi^T \Phi\), \( {\mb q}= {M}^{-1}\Phi^T \hat{\mb u},\)  then \(\widehat{\mb u}= \Phi \mb q\), and we have,
\begin{equation}\label{eq: Gal0'}
   M \ddot{\mb q}= \Phi^T \mb f(\Phi \mb q).
\end{equation}
For the linear problem, these equations become
\begin{equation}\label{eq: mass}
  M \ddot{\mb q}=-K\mb q,
\end{equation}
where
 \begin{equation}\label{eq: stiff}
    K= \Phi^T A \Phi. 
\end{equation}
The matrices \(M\) and \(K\) can be regarded as the mass and stiffness matrices in finite element approximations
of elasto-dynamics problems \cite{Hughes1987}.

\subsection{An error analysis}
The induced norm will be denoted by \(\| \cdot \|_0\). For a given time interval \([0,T]\), we also define,
\[ \| \mb u \|_{L^\infty(0,T;X)} = \sup_{t\in [0,T]} \|\mb u(t)\|_0,\]
and,
\[ \| \mb u \|_{L^2(0,T;X)} = \Big[ \int_0^T \|\mb u(t)\|_0^2 dt\Big]^\half.\]

In addition, we define the Ritz projection operator to the subspace, 
\begin{equation}
  \widehat{P}= \Phi (\Phi^T A \Phi)^{-1} \Phi^TA,
\end{equation}
with \(\widehat{Q}=I-\widehat{P}\) being the complementary projection. 

To compute the error, we let \(\mb e= \hat{\mb u} - \mb u= \theta + \eta\), in which,
\begin{equation}\label{eq: theta-eta}
  \begin{array}{rcl}
    \theta&=& \hat{\mb u} - \widehat{P}\mb u,\\
    \eta&=& \widehat{P}\mb u - \mb u= -\widehat{Q} \mb u.
  \end{array}
\end{equation}
For the initial condition, we will assume that 
\begin{equation}
 \hat{\mb u}(0)=P\mb u(0), \quad  \dot{\hat{\mb u}}(0)=P\dot{\mb u}(0).
\end{equation}

 The analysis is similar to the analysis of finite element methods
for second-order wave type of equations \cite{Baker76,BakDou76,Dup73,Geveci84,LaTh03}. The analysis presented here follows the approach in \cite{Baker76}, and provides the error estimate in the \(\|\cdot\|_0\) norm.

We consider the Galerkin method for the linear problem. In this case, 
the error equation for \(\theta\) in the variational form is given by,
\begin{equation}
  \big( \ddot{\theta}, \varphi \big) + a(\theta, \varphi)=   -\big( \ddot{\eta}, \varphi \big). 
\end{equation}
for any \(\varphi \in Y\).

We now re-write the equation as,
\[
  \frac{d}{dt} \big( \dot{\theta}, \varphi \big) -  \big( \dot{\theta}, \dot{\varphi} \big)
   = -a(\theta, \varphi) - \frac{d}{dt}\big( \dot{\eta}, \varphi \big)+\big( \dot{\eta}, \dot{\varphi}\big).\]
Next fix a \(\tau > 0\), we choose a particular test function as,
\begin{equation}
  \varphi(t)= \int_t^\tau \theta(s) ds,
\end{equation}
This special choice implies that
 \begin{equation}
 \varphi(\tau)=0, \quad \dot{\varphi}=-\theta.
\end{equation}
With this choice the equation above becomes,
\begin{equation}
\frac{d}{dt} \half   \big( {\theta}, \theta \big)
  =  -\frac{d}{dt} \big( \dot{\mb e}, \varphi \big) 
 + \frac{d}{dt} \half a(\varphi, \varphi) -\big( \dot{\eta}, \theta \big).
\end{equation}

Integrating the equation from \(0\) to \(\tau\), one gets,
\begin{equation}
 \begin{array}{ccl}
   \|\theta(\tau)\|^2_0 &=& \|\theta(0)\|^2_0
  - a\big(\varphi(0), \varphi(0)\big) \\
  &&\dsp - 2\big(\dot{\mb e}(0), \varphi(0)\big)- 2\int_0^\tau  \big( \dot{\eta}(s), \theta(s) \big) ds.
 \end{array}
\end{equation}
We will choose the initial condition so that \(\big(\dot{\mb e}(0), \varphi(0)\big)=0\). In addition, we notice that, \(a\big(\varphi(0), \varphi(0)\big) \ge 0\).  Therefore, using Cauchy-Schwartz inequality, we get for each \(\tau \in [0, T]\),
\begin{equation}
   \|\theta(\tau)\|^2_0 \le  \|\theta(0)\|^2_0 + \half \|\theta\|_{L^{\infty}(0,T;X)} + 2 T \|\dot{\eta}\|_{L^2(0,T;X)}^2. 
 \end{equation}
This leads to,
\begin{equation}
   \|\theta\|_{L^{\infty}(0,T;X)} \le \sqrt{2}\|\theta(0)\|_0  + 2\sqrt{ T} \|\dot{\eta}\|_{L^2(0,T;X)}. 
 \end{equation}

Combining with \eqref{eq: theta-eta}, we get,
\begin{Thm}
The error of the Galerkin method for the linear problem satisfies,
\begin{equation}
    \|\hat{\mb u} - \mb u\|_{L^{\infty}(0,T; X)} \le   C\Big( \|Q\mb u\|_{L^{\infty}(0,T;X)}   + \sqrt{ T} \|Q\dot{\mb u}\|_{L^2(0,T; X)} \Big). 
\end{equation}
\end{Thm}

\begin{rem}For perfect crystals, it is often possible to introduce norms that are similar to \(H^1(\Omega)\) and
\(H^2(\Omega)\) norms. In this case, one might be able to derive a more explicit bound for the interpolation  \cite{Lin03} and obtains an error estimate in \(H^1(\Omega)\) norm. In general, however, some symmetry may be lost near a lattice defect, it would be difficulty to define such norms.
\end{rem}

\begin{rem} There are two cases when such an estimate loses its values. The first case is when the time scale 
is  beyond \(\cal{O}(1)\) scale, \ie, the molecular times scale, which is about the scale of pico-seconds ($10^{-12}$s).  In fact, long time simulations of molecular systems is needed for most problems. Another problem
is when phonons are generated by local defects. Due to high frequency phonons, the term $\|Q\mb u\|$ might be large.  
\end{rem}

\subsection{A one-dimensional example}
To illustrate the idea of using the Galerkin projection as a reduction method, let's consider a one-dimensional 
model. 
\begin{equation}\label{eq: 1d-mod}
  m \ddot{x}_j=  V'(x_{j+1}-x_j) - V'(x_j-x_{j-1}). 
\end{equation}
We use the Lennard-Jones model for the potential energy $V$, given by,
\begin{equation}
 V(r) = r^{-12} - r^{-6}.
\end{equation}
The lattice spacing \(a_0\) and the mass are both normalized to unity.
We consider a chain with 1024 atoms. We will coarse-grain the first 512 atoms and keep the remaining 512 atoms.
This divides the system into two sub-domains. In the coarse-grained region, we choose one node out of every eight atoms, and we use usual piecewise linear basis function.   On the other hand, in the MD region, every atom is chosen as a node. We impose homogeneous boundary conditions at the ends. For the initial condition, we start with a left-moving wave packet centered in the domain on the right. To create such initial condition, we first
define,
\begin{equation}
 u^0(x)= \sum_{j=0}^{20} \cos(\xi_k x), \;   v^0(x)= \sum_{j=0}^{20} - \omega(\xi_k) \sin (\xi_k x),
\end{equation}
where $\omega(\xi)= 2 \sin \frac{\xi}{2}$ is the dispersion relation, and \(\xi_k= 0.5 + 0.02 k\) 
selects a few wave numbers. Then we choose the initial condition as follows,
\begin{equation}
  u_j(0) = F(j) u^{0}(j), \quad v_j(0) = F(j) v^{0}(j). 
\end{equation}
Here \(F(x) = 0.00025\exp\{-(x-640)^2/800\}\) is a Gaussian profile to confine the wave packet inside
the atomistic region.

\begin{figure}[htbp]
\begin{center}
\includegraphics[width=1.5in,height=5.5in,angle=-90]{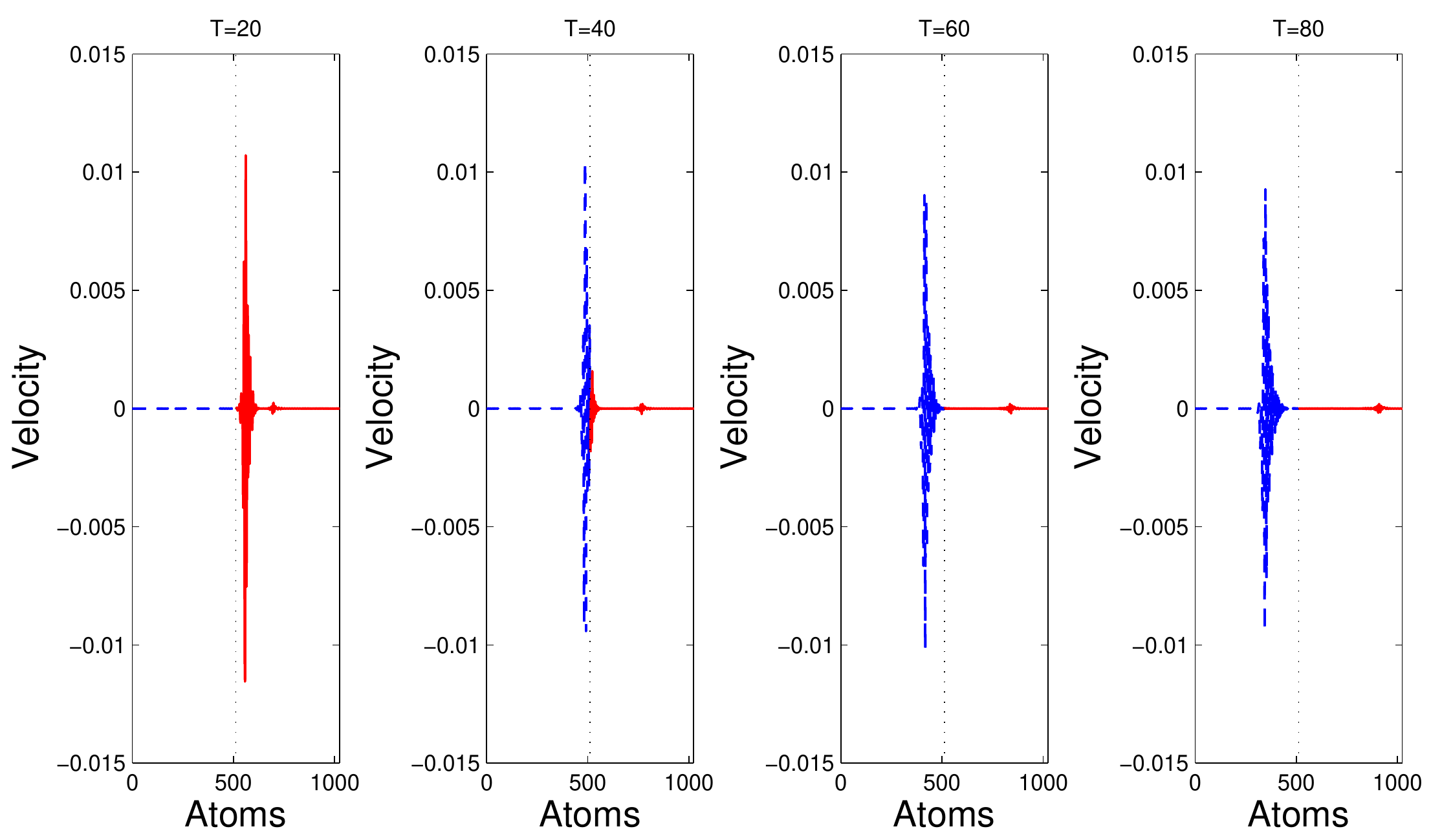}
\includegraphics[width=1.5in,height=5.5in,angle=-90]{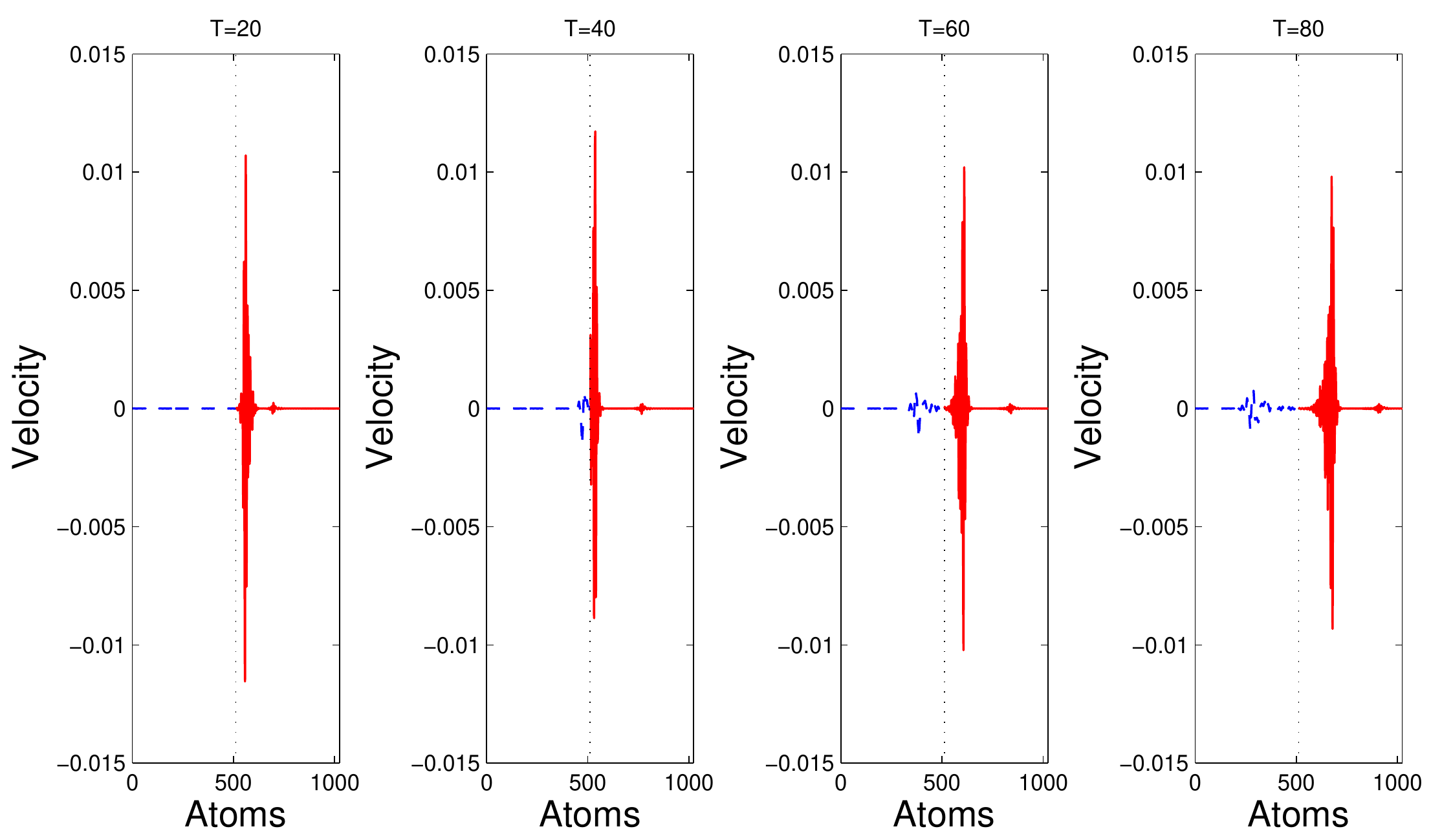}
\includegraphics[width=1.5in,height=5.5in,angle=-90]{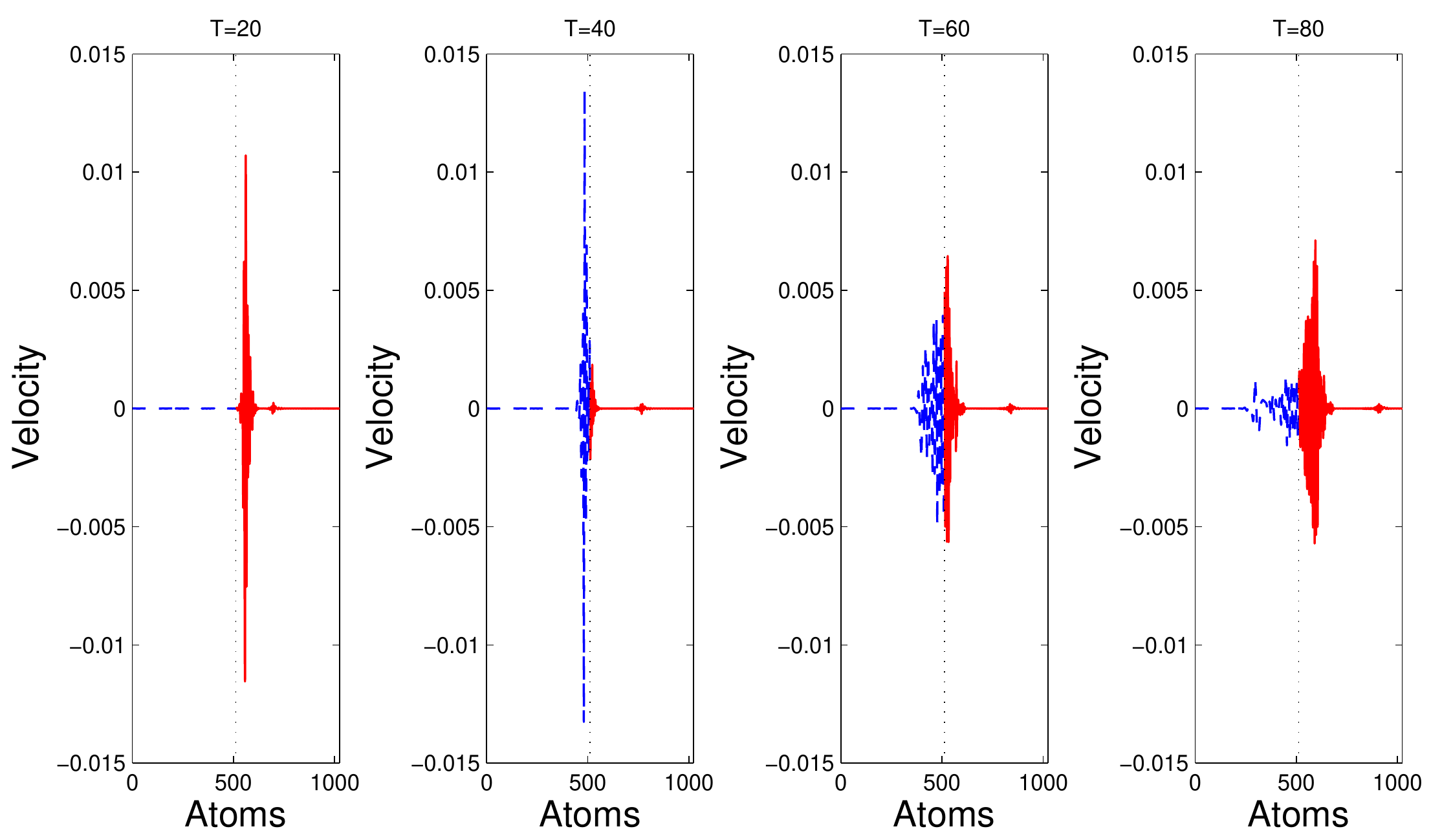}
\caption{Numerical results based on the one-dimensional model.}
\label{fig: chain012}
\end{center}
\caption{A numerical test with the one-dimensional model. From top to bottom: The results from a full atomistic simulation;  The results from the conventional Galerkin method; The results from the conventional Galerkin method with slowly varying mesh size. From left to right: Snapshots at T=20, 40, 60 and 80.} 
\end{figure}

The top figure in Figure \ref{fig: chain012} shows the result from a full atomistic simulation. We observe that the wave packet moves across the interface and enters the region on the left. For the conventional Galerkin method, however, the wave packet arrives at the interface, and then most of it is reflected back into the region on the right. 
We may attribute the reflection to the large ratio between the mesh size and the atomic spacing. Therefore, we conducted another test using the Galerkin method in which the mesh size changes linearly from $8a_0$ to $a_0$
as we move from the interface to the left boundary.  The figures on the third row in Figure \ref{fig: chain012} show the corresponding results. We see that although the mesh size varies very slowly, we still get significant reflections. 
Compared the case with uniform mesh size, the reflection occurred a bit later. 

\section{The Extended Galerkin Method}\label{sec: ext}

We have observed from the one-dimensional example that the conventional Galerkin method may give poor accuracy, especially when phonons are propagating out the atomistic region. In this section, we aim to improve the CG model \eqref{eq: Gal0'} derived from the conventional Galerkin method. 

\subsection{Galerkin approximation with memory}
We first derive an {\it exact} model for the coarse-grained variable \(\mb q\).  Let \(\mb y(t)= \Phi^T \mb u(t)\). We set our goal to finding an effective equation for \(\mb y\). To better illustrate 
 the ideas, we first consider the linear model. 
Taking time derivatives and split the right hand side according to the projection operators, we get,
\begin{equation}\label{eq: p-step1}
  \ddot{\mb y}= -\Phi^T A {P}\mb u - \Phi^T A  {Q} \mb u.
\end{equation}
The first term on the right hand side can be written as, 
\[-\Phi^T A {P}\mb u= - K M^{-1} \mb y.\]
Here \(M\) and \(K\) are the mass and stiffness matrices defined in \eqref{eq: mass} and \eqref{eq: stiff}.

Notice that by dropping the second term on the right hand side, and letting \(\mb q= M^{-1} \mb y\),  we obtain the conventional Galerkin approximation  \eqref{eq: Gal0} or \eqref{eq: Gal0'}. The example from the previous 
section suggests that a truncation at this step would be premature.   
To obtain a more accurate solution in time, we do not completely neglect the second term. Instead, we define, 
\begin{equation}
  \mb z= A {Q} \mb u.
\end{equation}
Again taking derivatives and split the right hand side in a similar way, we find,
\begin{equation}\label{eq: p-step2}
  \ddot{\mb z}= -AQ\mb z - AQAP\mb u = AQ\mb z - AQA\Phi M^{-1} \mb y.
\end{equation}
 To explicitly express \(\mb z\), we define 
fundamental solutions, \(C(t)\) and \(S(t)\), which satisfy the equations respectively,
\begin{equation}
  \ddot{C}= -AQ C, \; C(0)=I, \; \dot{C}(0)= 0,
\end{equation}
and
\begin{equation}
  \ddot{S}= -AQ S, \;S(0)=0, \;\dot{S}(0)= I.
\end{equation}
Clearly, these two functions both commute with \(AQ\). In addition, they are related as follows,
\begin{equation}\label{eq: C-S} 
\dot{S}(t)= C(t), \quad \dot{C}(t)=-S(t)AQ.
\end{equation}
Formally, the two functions can be written as,
\[ C(t)= \cos (AQ)^\half t, \; S(t)= (AQ)^{-\half}\sin (AQ)^\half t.\] 
Although the matrix \(AQ\) may not be invertible, the second expression can be defined by the corresponding
power series of the sinc function.

With these two functions, we get,
\begin{equation}
  \mb z(t)= C(t) \mb z(0) + S(t) \dot{\mb z}(0) - \int_0^t S(t-\tau) A{Q}A \Phi M^{-1} \mb y(\tau) d\tau.
\end{equation}

By collecting terms, we find that,
\begin{equation}
  \ddot{\mb y}= -\Phi^T A  {P}\mb u  +  \int_0^t \Phi^T S(t-\tau) A {Q}A \Phi M^{-1} \mb y(\tau) d\tau +\Phi^T  {\mb g}(t),
\end{equation}
where,
\begin{equation}
  \mb g(t)=  - \Big[ C(t) A {Q} \mb u(0) + S(t) A {Q} \dot{\mb u}(0)\Big].
\end{equation}

Finally, we let \( {\mb q}= {M}^{-1}\mb y\) and we obtain,
\begin{equation}\label{eq: q}
   {M} \ddot{  {\mb q} }(t) = -K  {\mb q}(t) + \int_0^t  {\Theta}(t-\tau) {\mb q}(\tau) d\tau +  \Phi^T {\mb g}(t).
\end{equation}
We let  \(\mb w= P \mb u=\Phi \mb q\). It satisfies an equation in a weak form,
\begin{equation}\label{eq: WF-exact}
  \big( \ddot{{\mb w}}, \mb v\big) + a( {\mb w}, \mb v) =   \int_0^t \big( B(t-\tau) {\mb w}(\tau), \mb v(t) \big) d\tau + \big({\mb g}, \mb v\big).
\end{equation}

Here the memory kernel functions are given by,
\begin{equation}
   {\Theta}(t)= \Phi^T B(t) \Phi, \quad B(t)= S(t) A {Q}A. 
\end{equation}
$\Theta $ is a symmetric, matrix-valued function.

The equation \eqref{eq: q} is an exact equation for \(\mb q\). We now make an approximation by neglecting \({\mb g}(t)\).
To justify this step, and to find more explicit expressions of the kernel function $\Theta$, we state the following lemma. 

\begin{Lem}
    Let \(\lambda \in \sigma(AQ)\) be an eigenvalue of \(AQ\), and \(\xi\) be the corresponding eigenvector.
    Then
    \begin{enumerate}
\item \(\lambda=0  \Longleftrightarrow \xi \in Y.\)
\item If \( \lambda\ne 0\), then \(\lambda\) is an eigenvalue of \(QAQ\), and the corresponding eigenvector is \(\eta=Q\xi\). 
\item If  \(\lambda\ne 0\) is an eigenvalue of \(QAQ\), and the corresponding eigenvector is \(\eta\in Y^\perp\), then \(\lambda\) is an eigenvalue of \(AQ\), and \(\xi= \eta + \lambda^{-1}PA \eta=  \lambda^{-1}A \eta\) is the corresponding eigenvector.
\item Let \(\lambda_i, i=1,2,\cdots,n\), be the eigenvalues of \(QAQ\), and \(\eta_i\) be the corresponding eigenvectors with unit length,
then
\[ QAQ= \sum_i \lambda_i \eta_i \eta_i^T, \quad Q= \sum_{\lambda_i\ne 0}\eta_i \eta_i^T, \;  \mathrm{and} \;AQ= \sum_{\lambda_i\ne 0} \lambda_i \xi_i \eta_i^T.\]
\item The force $\mb g(t)$ is expressed as,
\begin{equation}
 {\mb g}(t)=\Phi^T A\sum_{\lambda_i\in \sigma(QAQ)\backslash\{0\}} \cos \sqrt{\lambda_i} t\; \eta_i\eta^T_i \mb u(0)
 + \Phi^T A \sum_{\lambda_i\in \sigma(QAQ)\backslash\{0\}} \frac{\sin \sqrt{\lambda_i} t}{\sqrt{\lambda_i}} \; \eta_i\eta^T_i \dot{\mb u}(0).
\end{equation}

\item The history function \(\Theta(t)\) can be expressed as:
\[\Theta(t)= \Phi^TS(t)AQA\Phi= \sum_{\lambda_i\in \sigma(QAQ)\backslash\{0\}} \frac{\sin \sqrt{\lambda_i} t}{\sqrt{\lambda_i}} (\Phi^T A \eta_i)(\Phi^T A \eta_i)^T.\] 

\end{enumerate}

\end{Lem}
Notice that the zero eigen-modes of \(AQ\) do not appear in the above expressions because of the orthogonality conditions.

As a result, when the initial condition is close to the subspace \(Y\), \(\mb g(t)\) would be small. Therefore, we may neglect $\mb g(t)$ to obtain a closed system. In particular, when 
$\mb u(0)\perp Y$ and $\dot{\mb u}(0) \perp Y$, $\mb g(t) \equiv 0.$

This yields the following equation,
\begin{equation}\label{eq: Gal1}
   {M} \ddot{  {\mb q} }(t) = -K  {\mb q}(t) + \int_0^t  {\Theta}(t-\tau) {\mb q}(\tau) d\tau.
\end{equation}
The solution  \(\mb u\) is approximated by \( \hat{\mb u}  = \Phi \mb q\). The corresponding weak form is given by,
\begin{equation}\label{eq: WF-approx}
  \big( \ddot{\hat{\mb u}}, \mb v\big) + a( \hat{\mb u}, \mb v) =   \int_0^t \big( \Theta(t-\tau) \hat{\mb u}(\tau), \mb v(t) \big) d\tau,
\end{equation}
for each \(\mb v(t) \in Y\) with homogeneous boundary condition.  For nonlinear problems, one may replace the bilinear form \(a(\hat{\mb u},\mb v)\) by the nonlinear form \eqref{eq: nonlinear} and express the effective equation in a similar form as  \eqref{eq: Gal1}.

\begin{figure}[htbp]
\begin{center}
\includegraphics[scale=0.6]{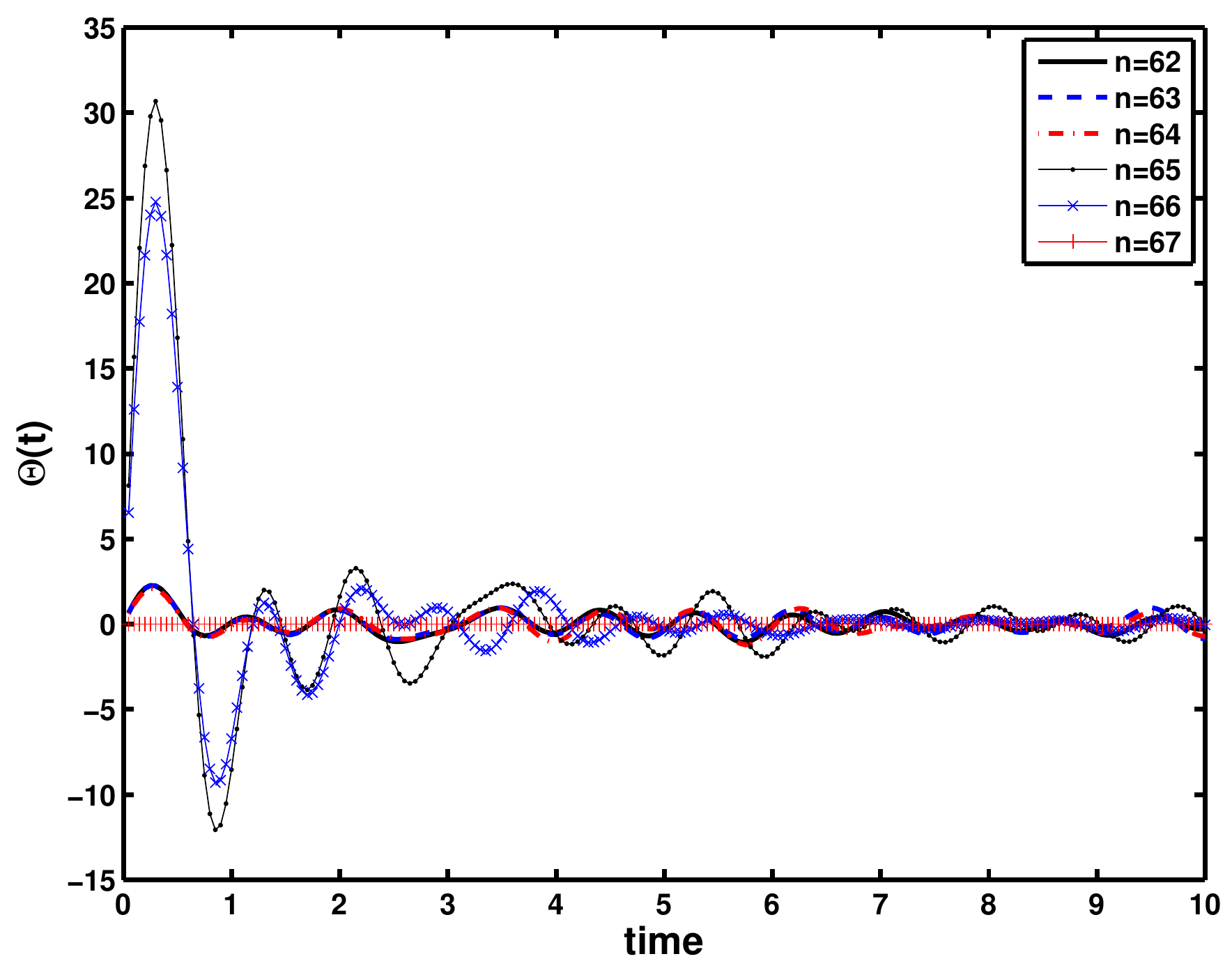}
\includegraphics[scale=0.85]{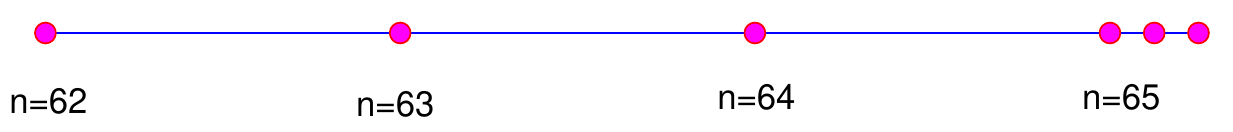}
\label{fig: theta}
\end{center}
\caption{The diagonal entries of the kernel function, $\Theta_{n,n}(t)$, for the one-dimensional example, with
$n=62, 63, \cdots, 67$. The bottom figure shows the corresponding nodes. $n=65$ is
the location of the interface.}
\end{figure}

We have computed the memory function $\Theta$  for the one-dimensional example described in the previous section. Figure \ref{fig: theta} shows the diagonal entries of the function
\(\Theta(t)\) for points near the interface. We can make the following observations: The kernel function
exhibits a peak near the interface, and away from the interface, the kernel function becomes significantly smaller. In particular, inside the MD region, it becomes zero. 

\subsection{Approximations of the GLE with direct expansion}

The GLE model provides a more accurate approximation by including a history-dependent term. In practice, however, the memory function has to be obtained prior to the computation. In general, direct calculation is quite difficult since the explicit expression involves solving \eqref{eq: C-S} which is almost as complicated as solving the original problem. Therefore the GLE often has to be approximated. One approach is to extend the system by introducing extended variables and embed the GLE into a system with no memory. This has been done for colored noise systems where the power spectrum is a rational function \cite{BaHaZh05}.  For our problem, we may rewrite the equation \eqref{eq: p-step1} into
\begin{equation}\label{eq: y}
  \ddot{\mb y}= -K M^{-1}\mb y - \Phi^T A  {Q} \mb u.
\end{equation}
The second term on the right hand side is not in the subspace $Y$. To proceed, we define,
\begin{equation}
  \bxi_0= \Phi^T A Q \mb u.
\end{equation}
Taking the second time-derivatives, and rearrange the terms on the right hand side, we can obtain,
\begin{equation}\label{eq: xi0}
  \ddot{\bxi}_0 = -\big[ K_2 - K_1 K^{-1}_0 K_1] \mb q + K_1 K_0^{-1} \bxi_0 - \Phi^T A^2 Q \mb u.
\end{equation}
Here, we have defined generalized stiffness matrices,
\begin{equation}
  K_\ell= \Phi^T A^\ell \Phi,
\end{equation}
for any \(\ell\ge 0\). Notice that \(K_0=M\), and its invertibility follows from the independence of the basis functions.

This equation contains a term that can be represented by neither \(\mb q\) nor \(\bxi_0\). So we define another additional variable,
\begin{equation}
 \bxi_1= \Phi^T A^2 Q \mb u. 
\end{equation}
This new variable satisfy the equation,
\begin{equation}
 \ddot{\bxi}_1= -\big[ K_3 - K_2 K^{-1}_0 K_1] \mb q + K_2 M_0^{-1} \bxi_0 - K_1 M_0^{-1} \bxi_1 - \Phi^T A Q A^2 Q \mb u.
\end{equation}

The last term in this equation is a new quantity. Conceivably, this process can be continued indefinitely, and a hierarchy can be built. This is similar to the BBGKY moment system for the Liouville equation. To obtain a closed finite system, truncation is often needed.  For instance, we may drop the last term in the equation above, and we find an extended system,
\begin{equation}
 \left\{
  \begin{array}{rcl}
      K_0 \ddot{\mb q}&=&  -K_1 \mb q - \bxi_0,\\
       \ddot{\bxi}_0 &=& -\big[ K_2 - K_1 K^{-1}_0 K_1] \mb q + K_1 K_0^{-1} \bxi_0 - \bxi_1,\\
        \ddot{\bxi}_1&=& -\big[ K_3 - K_2 K^{-1}_0 K_1] \mb q + K_2 K_0^{-1} \bxi_0 - K_1 K_0^{-1} \bxi_1.
  \end{array}
  \right.
\end{equation}   

We write this system in a more compact form,
 \begin{equation}
 \wt{M} \ddot{\wt{\mb q}} = - \wt{K} {\tilde{\mb q}},
\end{equation}
where,
\[
  {\wt{\mb q}} = 
  \left[
  \begin{array}{c}
  \mb q\\
  \bxi_0\\
  \bxi_1
  \end{array}
  \right], \quad
 \wt{M}=
  \left[
  \begin{array}{ccc}
  K_0 &0 &0\\
 0 &I &0\\
 0 &0& I 
  \end{array}
  \right], \quad
 \wt{K}=
  \left[
  \begin{array}{ccc}
  K_1 & I& 0\\
  K_2 - K_1 K^{-1}_0 K_1 &-K_1 K_0^{-1}  & I\\
  K_3 - K_2 K^{-1}_0 K_1&-K_2 K_0^{-1}  & K_1 K_0^{-1}
  \end{array}
  \right].\]
   
We further remark that,
\begin{rem}
 Such extension is not unique. For instance, one can instead define the additional variables as follows,
\[ \bxi_0= A Q \mb u,\; \bxi_1= (AQ)^2 \mb u, \;\bxi_2= (AQ)^3 \mb u, ...,\]
and the extended system reads,
\[
\begin{aligned}
\ddot{\bxi_0}&= - (K_2 - K_1 K_0^{-1} K_1 ) \mb y - \bxi_1,\\
\ddot{\bxi_1}&= -(K_3 - K_2 K_0^{-1} K_1 - K_1 K_0^{-1} K_2 + K_1 K_0^{-1} K_1 K_0^{-1} K_1 ) \mb y
- \bxi_2,\\
\cdots &
\end{aligned}
\]
\end{rem}
 
\begin{rem} 
  The additional variables may not have direct physical meanings as the system is further extended;
 \end{rem}
 \begin{rem}
 The variational structure is lost, \ie, there is no weak form associated with the truncated system, and the extended system may not be well-posed.
\end{rem}

To demonstrate the potential ill-posedness of the extended system, let's consider the one-dimensional chain model with nearest neighbors. We choose a chain of 5 atoms with 
indices from 0 to 4, and we fix the first and the last atom to their equilibrium position, \ie, $u_0=0, u_4=0$. 
With some normalization, we choose the matrix \(A\) as follows,
\[
A=\left[
\begin{array}{cccc}
 2 & -1 & 0 & 0\\
 -1 & 2 & -1 & 0 \\
 0 & -1 & 2 & -1 \\
 0 & 0 & -1 & 2
\end{array}
\right].\]
This matrix is positive definite. Therefore the trivial solution of the linearized MD model is stable.
We now choose one CG variable with nodal value defined at the third atom, and we set,
\[
\Phi^T=[ \half,\; 1,\; \half],\]
which corresponds to a linear interpolation.  The stiffness matrices can be easily computed. They are given by,
\[K_0=\frac{3}{\;2 \;},\; K_1=K_2=1.\]
We consider the extension to \(\bxi_0\), by dropping \(\bxi_1\) from the equation \eqref{eq: xi0}. In this case, we found that the generalized eigenvalues for \((\tilde{K}, \tilde{M})\) to be $\pm \sqrt{6}/3$. The negative eigenvalue
indicates instability.

\subsection{An extended Galerkin's approximation} 

From the example at the end of last section, one finds that such direct expansion is not reliable since it may introduce unstable models.  We propose a new extension method which maintain the Galerkin formulation. This extension is done by expanding the approximation space.  This would eliminate the problem of instability. One idea would be to extending the subspace by introducing bubble elements \cite{ArBrFo84,BaBrFr93}. In particular, this idea has been applied to multi scale problems \cite{Hughes99}. But here we will follow a different path.

For the purpose of the extension, we let \(Y_0=\mathrm{Range}(\Phi)\). We let \(Y_1=\mathrm{Range}(QA\Phi)\), and 
\(Y= Y_0 \oplus Y_1\). Then the Galerkin projection of the linear problem to \(Y\) yields 
\begin{equation}\label{eq: eG0}
\left\{
  \begin{aligned}
    K_0 \ddot{\mb q}&= - K_1 \mb q -  \wt{K}_2\mb \bxi_0\\
    \wt{K}_2 \ddot{\bxi}_0 &= -\wt{K}_2\mb q - \wt{K}_3 \bxi_0.
  \end{aligned}
  \right.
\end{equation}
Here, the approximate solution is expressed as, \(\hat{\mb u}=\Phi \mb q + Q A \Phi \bxi_0\). The matrices
\(\wt{K}_2\) and \(\wt{K}_3\) are given by,
\begin{align}
\wt{K}_2&= \Phi^T AQA\Phi = K_2 - K_1K_0^{-1}K_1, \\
\wt{K}_3&=\Phi^T AQAQA\Phi=K_3 - K_2 K_0^{-1} K_1 - K_1 K_0^{-1} K_2 + K_1 K_0^{-1} K_1 K_0^{-1} K_1 . 
\end{align}

We now go further to extend the approximation space. For example, we can introduce another projection,
\begin{equation} 
P_1 = QA\Phi (\Phi^TAQA\Phi)^{-1} \Phi^T A Q, \quad Q_1=I - P_1.
\end{equation}
We let \(Y_2= \mathrm{Range}(QQ_1A^2\Phi)\). One can easily verify that \(Y_2 \perp Y_0\) and \(Y_2\perp Y_1\).
We set \(Y=Y_0 \oplus Y_1\oplus Y_2\), and seek the solution in the form of,
 \begin{equation}
  \widehat{\mb u}= \Phi \mb q(t) + QA\Phi \bxi_0(t) + QQ_1A^2\Phi \bxi_1(t).
\end{equation}
Using the variational approach, we find,
\begin{equation}\label{eq: eG1}
\left\{
  \begin{aligned}
    K_0 \ddot{\mb q} &= -K_1 \mb q - \wt{K}_2 \bxi_0, \\
    \wt{K}_2 \ddot{\bxi}_0 &= -\wt{K}_2 \mb q - \wt{K}_3 \bxi_0 - \wt{K}_4\bxi_1, \\
    \wt{K}_4 \ddot{\bxi}_1 &= - \wt{K}_4 \bxi_0 - \wt{K}_5 \bxi_1. 
  \end{aligned}
  \right.
\end{equation}
Here the matrices are given by,
\[
\begin{aligned}
  \wt{K}_4&= \Phi^TAQAQQ_1 AQA\Phi, \\
  \wt{K}_5&= \Phi^TAQAQ_1QAQQ_1AQA\Phi.
\end{aligned}
\]
For nonlinear problems, the Galerkin projection will be applied to \eqref{eq: md0}. The variational procedure 
yields,
\[
\left\{
\begin{aligned}
 K_0 \ddot{\mb q} &= \Phi^T \mb f(\hat{\mb u}),\\
  \wt{K}_2 \ddot{\bxi}_0 &= \Phi^T A Q \mb f(\hat{\mb u}),\\ 
   \wt{K}_4 \ddot{\bxi}_1 &= \Phi^T A^2 Q_1 Q\mb f(\hat{\mb u}),
\end{aligned}
\right.
\] 
\smallskip

For the linearized model, the coarse-grained model derived from this extended Galerkin method can be viewed as an approximation
of the generalized Langevin model  \eqref{eq: Gal1}. To see this, we assume that \(\mb u(0) \in Y\) and
\(\dot{\mb u}(0) \in Y\), and in light of this, we set \(\bxi_0(0)=0\). Applying Laplace transform to the second equation in \eqref{eq: eG0}, we found that, 
\[
 -\wt{K}_2\bxi_0(t)= \int_0^t  \Theta_1(t-\tau) \mb q(\tau)d\tau,\]
 which can be substitute into the first equation, yielding,
\begin{equation}
  M\ddot{\mb q}= - K \mb q  + \int_0^t \Theta_1(t-s) \mb q(s)ds.
\end{equation}
Here the function \(\Theta_1\) has Laplace transform,
\begin{equation}
  \wt{\Theta}_1(\lambda)= \wt{K}_2 \bigl(\lambda^2 \wt{K}_2 + \wt{K}_3\bigr)^{-1} \wt{K}_2.
\end{equation}

Similarly, \eqref{eq: eG1} can be turned into the same form with the kernel function's Laplace
transform given by,
\begin{equation}
\wt{\Theta}_2(s)= \wt{K}_2\Big[ (s^2 \wt{K}_2 + \wt{K}_3) - \wt{K}_4 \big(s^2 \wt{K}_4 + \wt{K}_5\big)^{-1} \wt{K}_4
\Big] \wt{K}_2.
\end{equation}
This provides a rational approximation of the Laplace transform of the kernel function.

We briefly explain the motivation for using such subspaces for the extension procedure. Notice that when the initial
data are in the subspace \(Y\), the solution of the harmonic model \eqref{eq: md1},
is in the space, $K= \text{span}\{ A^\ell \Phi, \;\ell \ge 0\}$. This can be seen by writing down the solutions in the form of
sine and cosine functions, and then expanding the solution into power series. Therefore, we may consider the subspace \(Y= {\text{Range}}(\Phi) + {\text{Range}}(A\Phi)=Y_0 \oplus Y_1\). Similarly, one can directly verify that
\(Y= {\text{Range}}(\Phi) + {\text{Range}}(A\Phi) + {\text{Range}}(A^2\Phi)=Y_0 \oplus Y_1 \oplus Y_2.\)
To continue with the expansion, one may consider the sequence of subspaces $K_{L}(A, \Phi)=\text{span}\{
\Phi, A\Phi, \cdots, A^L \Phi\}$, which is a Krylov type of subspace. The Cayley-Hamilton theory implies that
$K_\ell \subseteq K_n$ for \(\ell \ge n\); $n$ is the dimension of the problem. If \(Y_0\) is not contained in the sum of the eigenspaces associated with
a subset of the eigenvalues, then $K_n = X$, and the exact solution should be recovered if the subspace $Y$ is expanded to $K_n$.

\subsection{The one-dimensional chain problem}

\begin{figure}[htbp]
\begin{center}
\includegraphics[width=1.5in,height=5.5in,angle=-90]{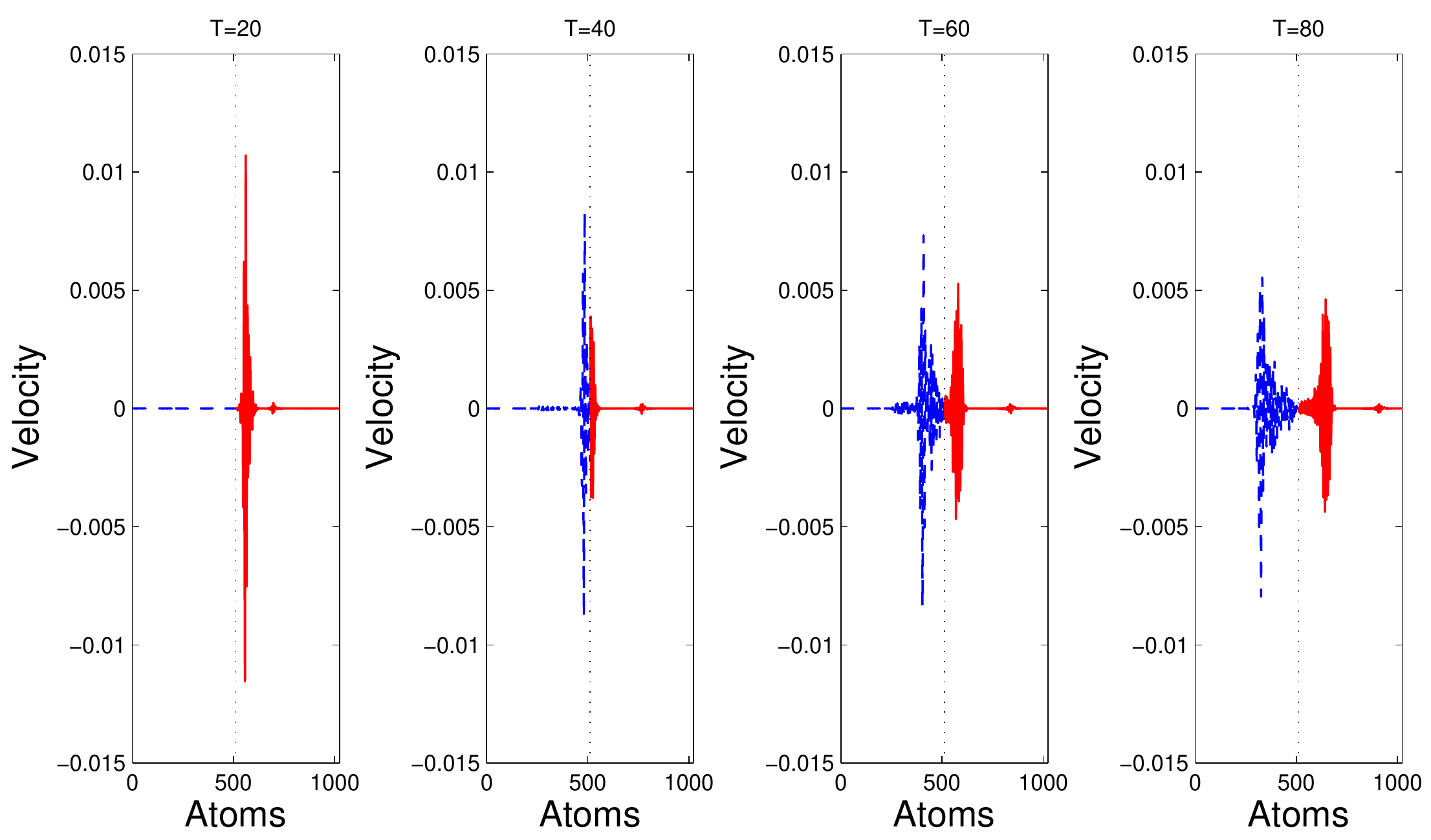}
\includegraphics[width=1.5in,height=5.5in,angle=-90]{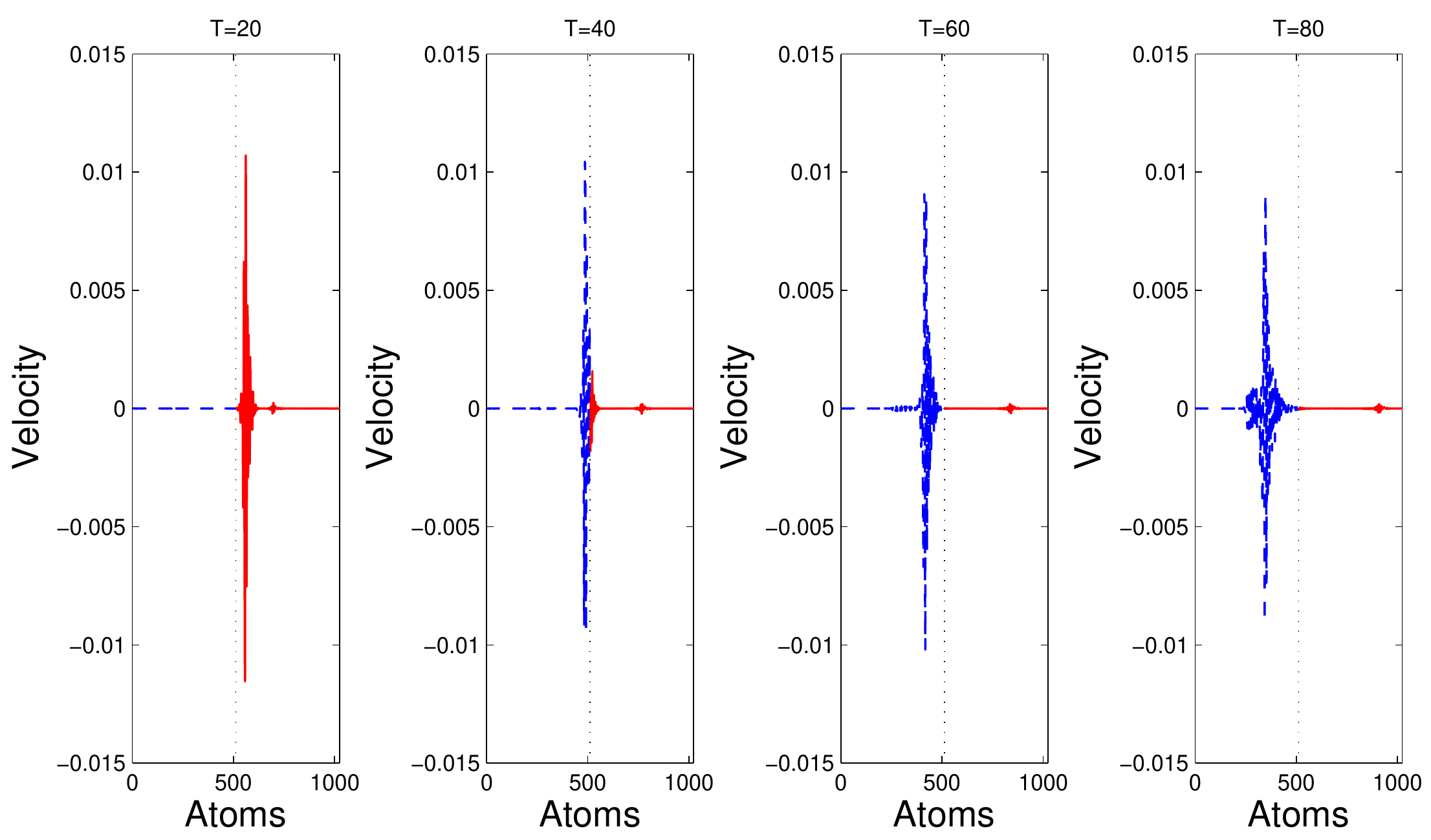}
\label{fig: chain3}
\end{center}
\caption{The solution of the 1d chain model using the extended Galerkin method. Top: $L=5$; Bottom: $L=10$.}
\end{figure}

We now turn again to the one-dimensional example. As in the conventional Galerkin method, we partition the system into a CG region where the mesh size $H=8a_0$.  For the extended Galerkin method, we need to 
select additional nodes near the interface. We choose 20
nodal points in the CG region, and 5 atoms from the atomistic region. The results are included in Fig. \ref{fig: chain3}. For the choice $L=5$, some reflections are observed. When $L$ is increase to 10, almost no reflection is observed. 

\section{Implementation Issues}\label{eq: imp}

In the formulation of the extended Galerkin method, we have assumed that all the quantities, \eg, the projection matrices and the stiffness matrices, can be computed  {\it exactly}. Clearly, computing these quantities exactly would require an effort that is as much as solving the original problem.  

\subsection{Construction of the approximation spaces}
We first discuss the construction of the subspace to which the Galerkin method is applied. 
Several steps are taken to simplify the procedure. 

First we discuss the subspace \(Y_0\). In light of the lattice defects in the system, where the molecular 
dynamics model should be kept, it is convenient to partition the system into a coarse-grained region, and an atomistic region. In the CG region, we use a conventional finite element representation of the displacement field, 
\eg,  piecewise linear functions. Inside the atomistic region, we pick the trivial basis function, \(\varphi(\mb x)= \delta_{nk}\) at the point \(\mb y^0_n\). This function is one at this point, and zero at all other points.  The union of such point will be donated by \(\Omega_{\text{at}}\). From the previous section, we see that
the memory function $\Theta(t)$ is small away from the interface. Therefore, we will only select the basis functions
with nodal points near the interface to expand the subspace. For each selected basis function, we let 
\begin{equation}\label{eq: krylov}
K_L(A,\varphi)= \text{span}\left\{ \varphi, A\varphi, \cdots, A^L \varphi\right\}. 
\end{equation}
The whole extended subspace is constructed by combing  $K_L(A;\varphi_j), j\in J$, where
$J$ is the set of selected basis functions. For this purpose, 
We need to include some more atoms  around \(\Omega_{\text{at}}\). The new set is denoted by \(\overline{\Omega}_{\text{at}}\). This will be explained next.

 We now discuss how to generate \(K_L(A,\varphi)\). We need to compute a vector in the form of
  $A\psi$ where \(\psi\) is the discrete function whose support lies in \(\overline{\Omega}_{\text{at}}\). The vector will be approximated by a finite difference formula,
\[  A\psi \approx -\frac{\mb f(\mb y + \veps \psi) - \mb f(\mb y )}{\veps}.\] Here \(\mb y\) is the  reference position.  For most empirical potentials, the interaction has finite range, and the shape function \(\psi\) has support in \(\overline{\Omega}_{\text{at}}\), $A\psi$ also has support close to $A\psi$. In our implementation,
we will choose \(\overline{\Omega}_{\text{at}}\) to contain the support of all \(A^\ell \varphi_j, \ell \le L, j\in J\).
 
Finally, the new set of functions generated from the previous step may not be independent. One implication is that the mass matrix  may be singular. In this case, we will use an orthogonalization procedure similar to the Anoldi's method \cite{Saad}. However, we only make the additional basis function orthogonal to $Y_0$.   This way, only the position of the atoms in \(\overline{\Omega}_{\text{at}}\)  are needed in the process.
We will denote the additional set of basis functions as, 
\begin{equation}
 \Psi= \left[ \psi_1, \psi_2, \cdots, \psi_k\right]. 
\end{equation}
 This is a set of orthonormal basis vectors. In addition, \(\Psi^T \Phi=0\).  
With this notation, we have \( Y= Y_0 \oplus \text{Range}(\Psi)\), and the Galerkin projection yields,
\begin{equation}\label{eq: cg}
  \begin{aligned}
  M \ddot{{\mb q}} &= \Phi^T \mb f(\hat{\mb u}),\\
    \ddot{\bxi} & = \Psi^T  \mb f(\hat{\mb u}),
  \end{aligned}
\end{equation}
where $\hat{\mb u}= \Phi \mb q + \Psi \bxi.$

\subsection{Summation and quadrature rules}
The coarse-grained MD model is an autonomous, second order ODEs. We use the Verlet's method to integrate the system in time. One advantage of this is that at each step of the time integration, the forces only have to be sampled once.  

We now discuss the calculation of the mass matrix and the averaged forces. 
As an example, we consider 
the case where the CG region has been divided in triangles (2D) or tetrahedrons (3D), with basis functions \(\varphi_n\) given by standard piecewise linear functions. 
For the mass matrix, we need to compute 
\begin{equation}\label{eq: mmn}
M_{m,n}= \sum_{i} \varphi_m(\mb y_i) \varphi_n(\mb y_i), 
\end{equation}
and for the forces, we need to compute,
\begin{equation}\label{eq: Fn}
F_n=\bigl(\varphi_n(\mb y), \mb f(\mb y)\bigr),
\end{equation}
 where \(\varphi_n\) is one of the basis functions.
 
For the nodal functions in $\Phi$, they are divided into two groups, the ones in  \(\overline{\Omega}_{\text{at}}\), which includes the basis functions in $\Psi$, and the ones in the coarse-grain region. For the formulas \eqref{eq: mmn} and \eqref{eq: Fn}, if the basis functions have support in \(\overline{\Omega}_{\text{at}}\), we will evaluate
the expression exactly. Otherwise, we will approximate the summation by integrals. This will be discussed next.

We first consider the computation of the mass matrix. Let \(m\) and \(n\) be two vertices of a triangle \(\mathcal{T}\) (or a tetrahedron in 3D). 
The sum in \eqref{eq: mmn} is over the intersection of the support of \(\varphi_m\) and \(\varphi_n\).  It suffices to consider the sum over a triangle, which can be approximated by an integral, as was done in \cite{WaLi03},
\begin{equation}
\sum_{i\in \cal{T}} \varphi_m(\mb y_i) \varphi_n(\mb y_i)
\approx \frac{1}{\cal{V}_0}\int_\cal{T} \varphi_m(\mb y) \varphi_n(\mb y) d\mb y.
\end{equation}
Here \(\cal{V}_0\) is a volume of the unit cell. This integral then can be further approximated by standard quadrature formulas.   

Next we consider the sampling of forces, given by the inner product, 
\[
 F_n=\bigl(\varphi_n(\mb y), \mb f(\mb y)\bigr).
\]
Quadrature approximation of such average has been discussed in \cite{WaLi03,GuZh10}. 
Here we will derive an expression that can be directly linked to the continuum stress. Let \(\mb f_i\) be the force exerted on the \(i\)th atom, and we assume that \(\mb f_i\) can
be decomposed as follows,
\begin{equation}\label{eq: fsum}
\mb f_i = \sum_{j\ne i} \mb f_{ij}, \quad \mb f_{ij}= - \mb f_{ji}.
\end{equation}
To our knowledge, all the empirical potentials for crystalline solids admit natural decompositions in this form.

Next, within the support of $\varphi_n$ denoted by \(\Omega_n\), we have,
\begin{align}
&\bigl(\varphi_n(\mb y), \mb f(\mb y)\bigr)\\
&= \sum_{i\in \Omega_n} \mb f_i \varphi_n(\mb y_i) \\
&= \sum_{i\in \Omega_n} \sum_{j\ne i}  \mb f_{ij} \varphi_n(\mb y_i)\\
&= \sum_{i\in \Omega_n} \sum_{\stackrel{j \in \Omega_n}{j\ne i}}\mb f_{ij} \varphi_n(\mb y_i) 
   +	\sum_{i\in \Omega_n} \sum_{\stackrel{j \notin \Omega_n}{j\ne i}} \mb f_{ij}\varphi_n(\mb y_i) \\
&= \half \sum_{i\in \Omega_n} \sum_{\stackrel{j \in \Omega_n}{j\ne i}}  \mb f_{ij}\bigl[\varphi_n(\mb y_i) - \varphi_n(\mb y_j)\bigr]
   +	\sum_{i\in \Omega_n} \sum_{{j \notin \Omega_n}}  \mb f_{ij}\varphi_n(\mb y_i)\\
&\approx -\sum_{i\in \Omega_n}  {\boldsymbol\sigma}_i \nabla \varphi_n(\mb y_i) \cal{V}_0 
+ \sum_{i\in \Omega_n} \sum_{{j \notin \Omega_n}}
 \mb f_{ij}\Bigl[ \varphi_n(\mb y_i) + \half \bigl(\mb y_j - \mb y_i) \cdot\nabla \varphi_n(\mb y_i)  \Bigr].   
\end{align}
Here we have defined a local stress,
\begin{equation}
{\boldsymbol \sigma}_i \overset{\rm def}{=} -\frac{1}{2\cal{V}_0}\sum_{j\ne i}  \mb f_{ij}\otimes(\mb y_i - \mb y_j).
\end{equation}
This has been motivated by the expression of virial stress \cite{zhou2003new,zimmerman2004calculation}. 

The last term, 
\(
\Big[\varphi_n(\mb y_i) + \half \bigl(\mb y_j - \mb y_i) \cdot \nabla \varphi_n(\mb y_i)  \Bigr], 
\) can be regarded as an extrapolation of the basis function to the boundary of \(\Omega_n\), where
the function vanishes. Therefore this term will be neglected.

With these approximations, the sum is reduced to,
\begin{equation}
\bigl(\varphi_n(\mb y), \mb f(\mb y)\bigr)\approx -\sum_{i\in \Omega_n} {\boldsymbol\sigma}_i\nabla \varphi_n(\mb y_i) \cal{V}_0. 
\end{equation}

In the implementation, this will be approximated by an integral, and further approximated by
mid-point quadrature formula in each triangle,
\begin{equation}
\sum_{i\in \Omega_n} {\boldsymbol\sigma}_i\nabla \varphi_n(\mb y_i) \cal{V}_0
\approx \int_{\Omega_n}{\boldsymbol\sigma}(\mb y) \nabla\varphi_n(\mb y) d\mb y.
\end{equation}

At this point, one might notice that with these approximations, the procedure in the CG region is the same as solving the elasto-dynamics equation,
\begin{equation}
 \rho_0 {\mb u}_{tt} = \nabla \cdot \sigma(\nabla \mb u),
\end{equation}
with the stress-strain relation given by the Cauchy-Born rule \cite{CB84}.

\section{Another Numerical example}\label{sec: num}

Here we consider a test problem for a dislocation dipole in aluminum. The atomic potential used here is the EAM potential \cite{ErAd94} with lattice constant $4.032$\AA. The system under consideration is a $3D$ rectangular sample, with the three orthogonal axes being along the $\langle110\rangle$, $\langle001\rangle$ and $\langle1\bar{1}0\rangle$ directions respectively. The system is periodic in the third direction with period equal to 8 times the atomic spacing.  Therefore, we will treat it as a two-dimensional problem in the coarse-graining procedure.

The dimension of the entire system is \( 45.62\; nm \times 48.38\; nm\), and  the full system contains 384,000 atoms.  Two dislocations are introduced with core position at \((\pm 30,0)\)\AA $\,$and opposite Burgers  vectors,
given by \(\mb b=\pm \half[1 1 0]\). The initial position of the atoms, denoted by \(\mb z_i\), is computed from superimposed analytical solutions of the anisotropic elasticity problems around edge dislocations \cite{Stroh58,Ting1997}. Meanwhile, a pure shear strain with \(e_{21}=0.02\) is applied at the remote boundary. This shear strain is also added to the atomic positions, and the new position will be set as the reference position of the atoms. Namely, we set  \(\mb y_i= \mb z_i + E \mb z_i\) with \(E\) being the strain tensor. The initial velocity is zero. For the time integration, we choose the step size \(\Delta t = 0.052880\) (pico-second).

\begin{figure}[htp]
\begin{center}
\includegraphics[scale=0.4]{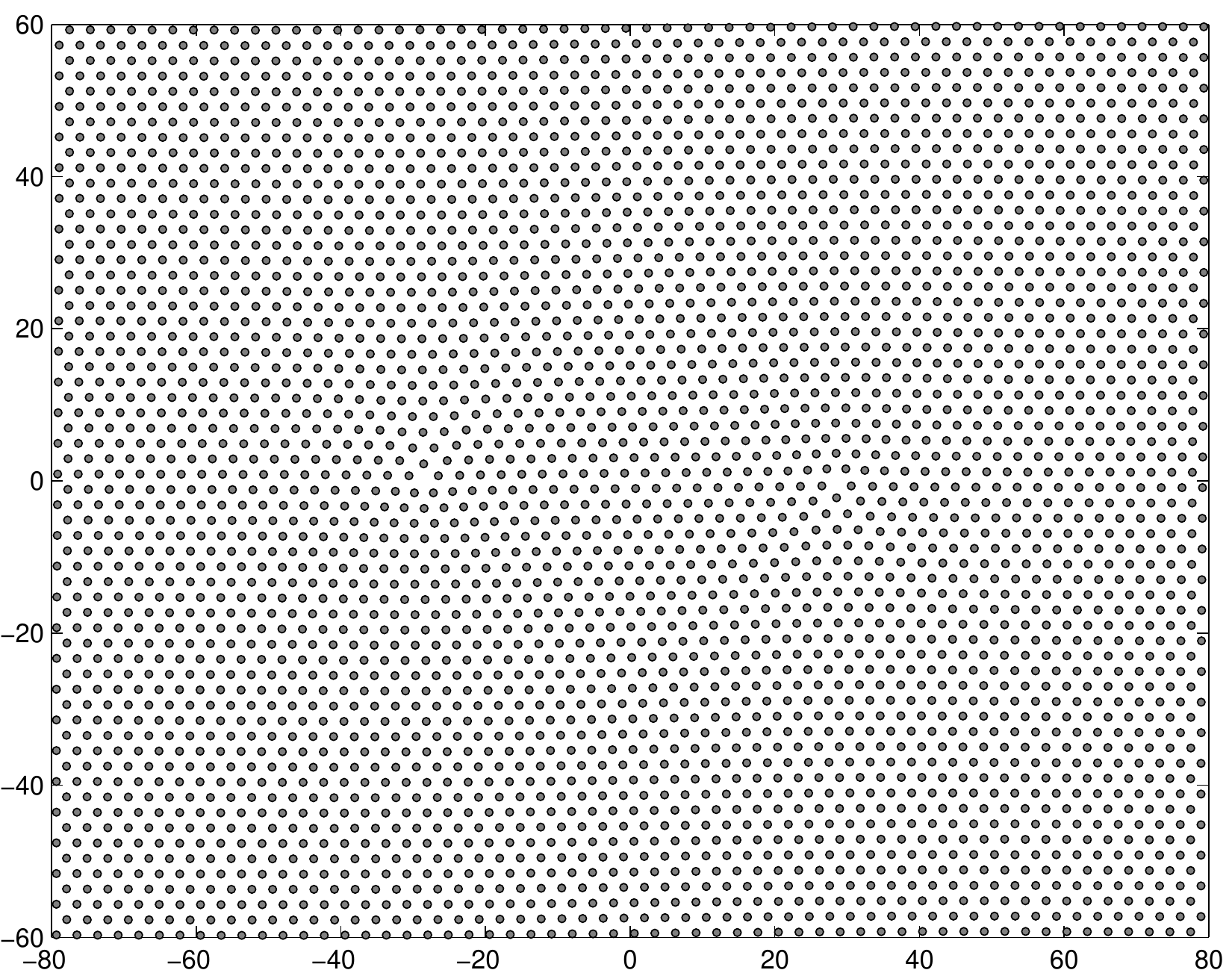}\\
\includegraphics[scale=0.4]{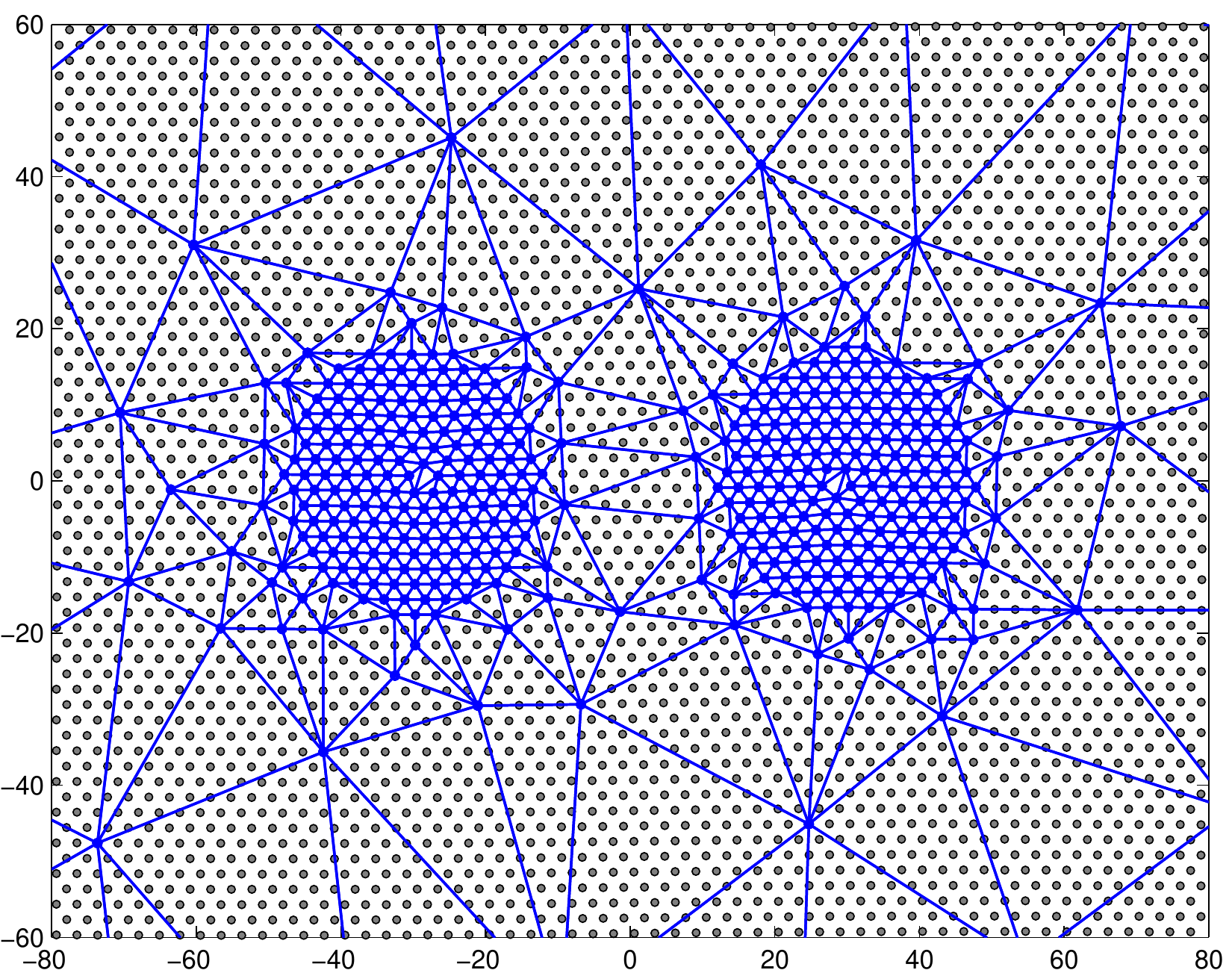}\\
\includegraphics[scale=0.4]{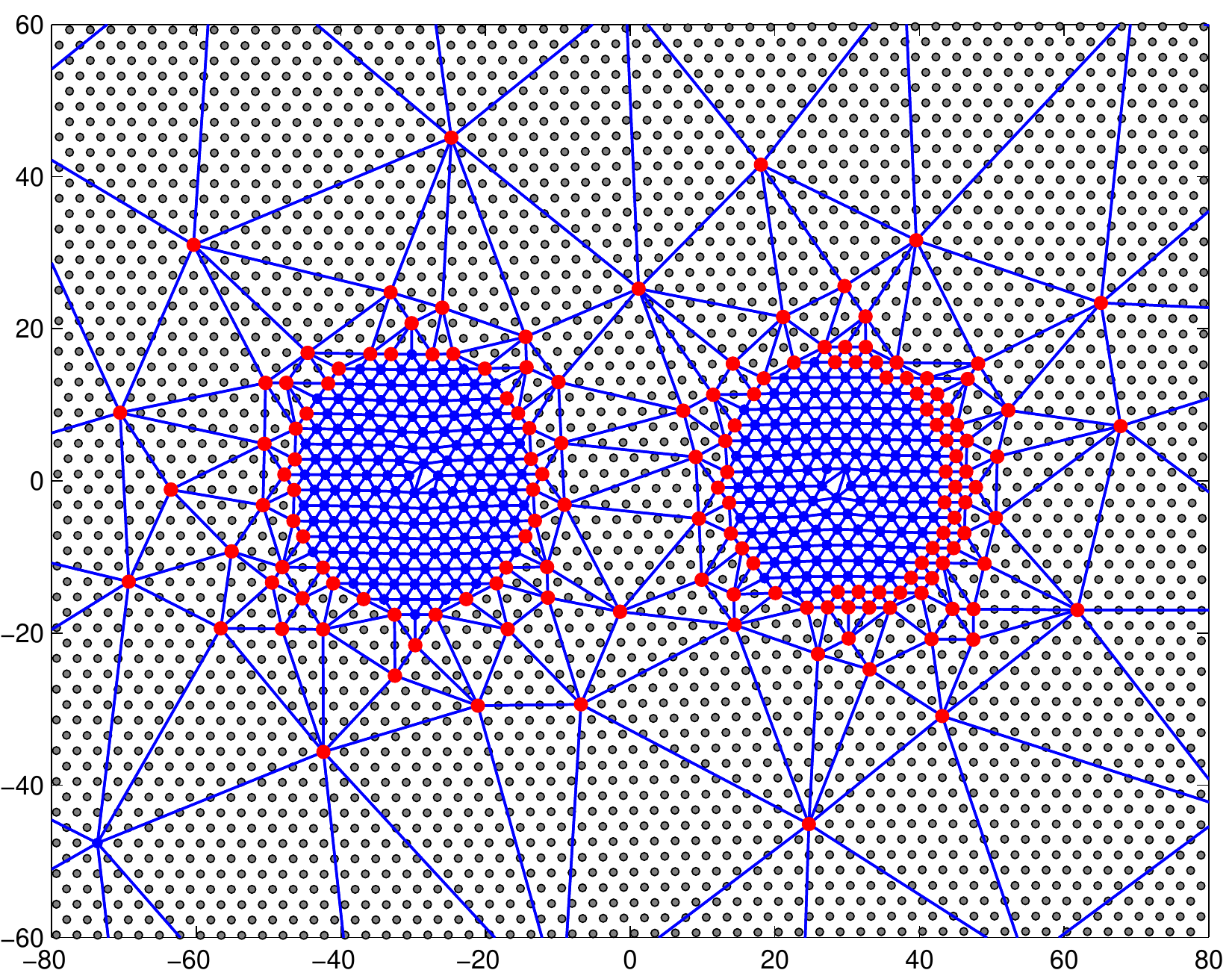}
\caption{Finite element nodes. Top: All the atoms in the system; Middle: the selected model point
for a conventional Galerkin method; Bottom: Highlighting the repeated nodes at the interface.}
\label{fig: mesh}
\end{center}
\end{figure}

In figure \eqref{fig: mesh}, one can see the projected position of the atoms as well as the finite element mesh 
around the two dislocation cores. In addition, in the bottom figure, we have indicated the selected nodes, where
the subspace is extended. These nodes will be regarded as {\it repeated nodes} since additional degrees of freedom are introduced at those points.

\begin{figure}[htb]
\begin{center}
\includegraphics[width=5in,height=2in]{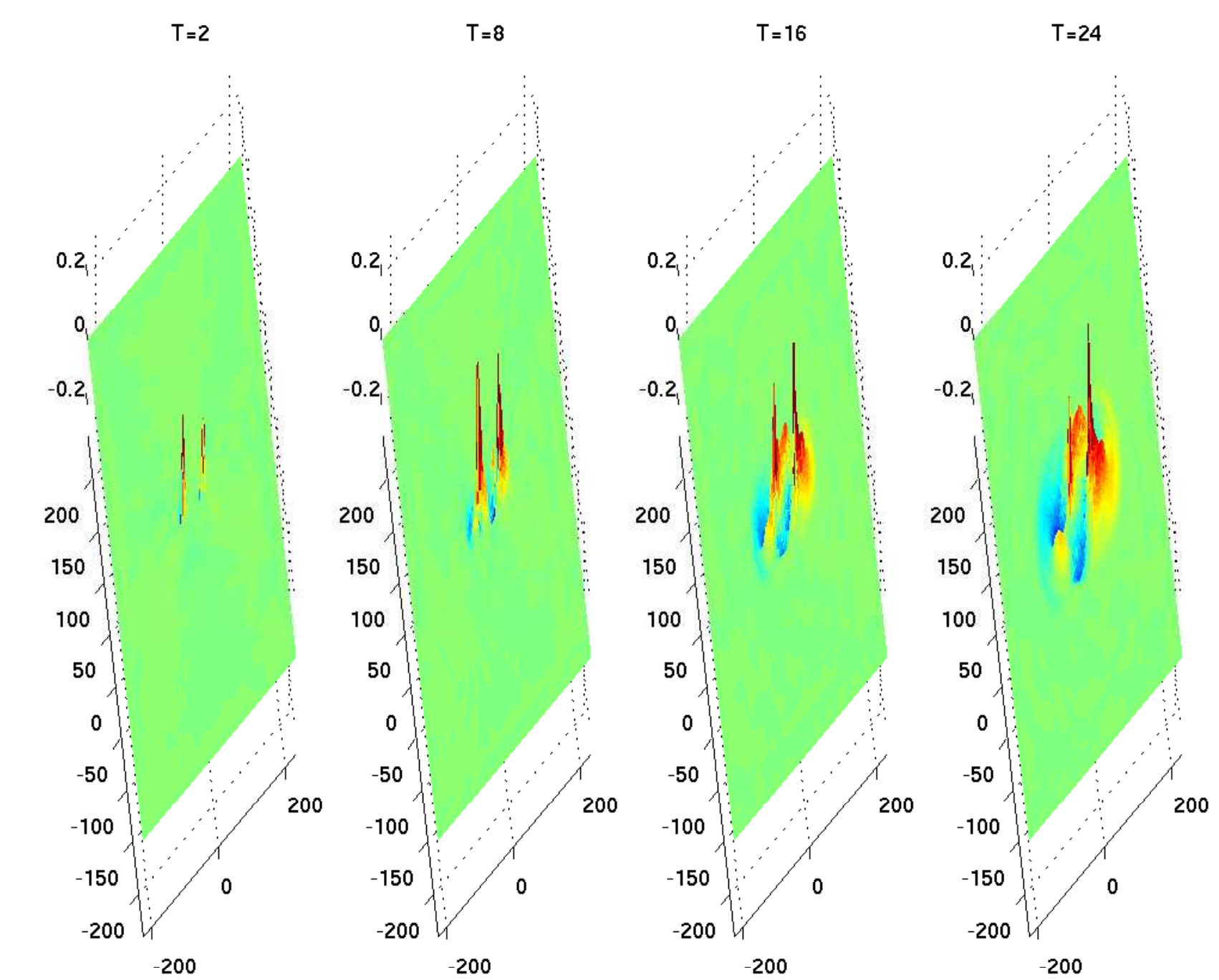}
\end{center}
\caption{The snapshots of the displacement \(u_1\) computed from the full molecular dynamics model. }
\label{fig: dd00}
\end{figure}

We first show the results from the full molecular dynamics model in Fig. \ref{fig: dd00}. Since the initial
position is determined from the corresponding continuum model, we expect that away from the dislocation core,
the atoms are almost at equilibrium at the beginning. Near the dislocation cores, the analytical solutions 
are not accurate due to the large deformation and loss of symmetry. This is confirmed by the
numerical results: The change of displacement
starts from the dislocation cores and then propagates into the surrounding region.  

\begin{figure}[htb]
\begin{center}
\includegraphics[width=1.25in,height=3.5in]{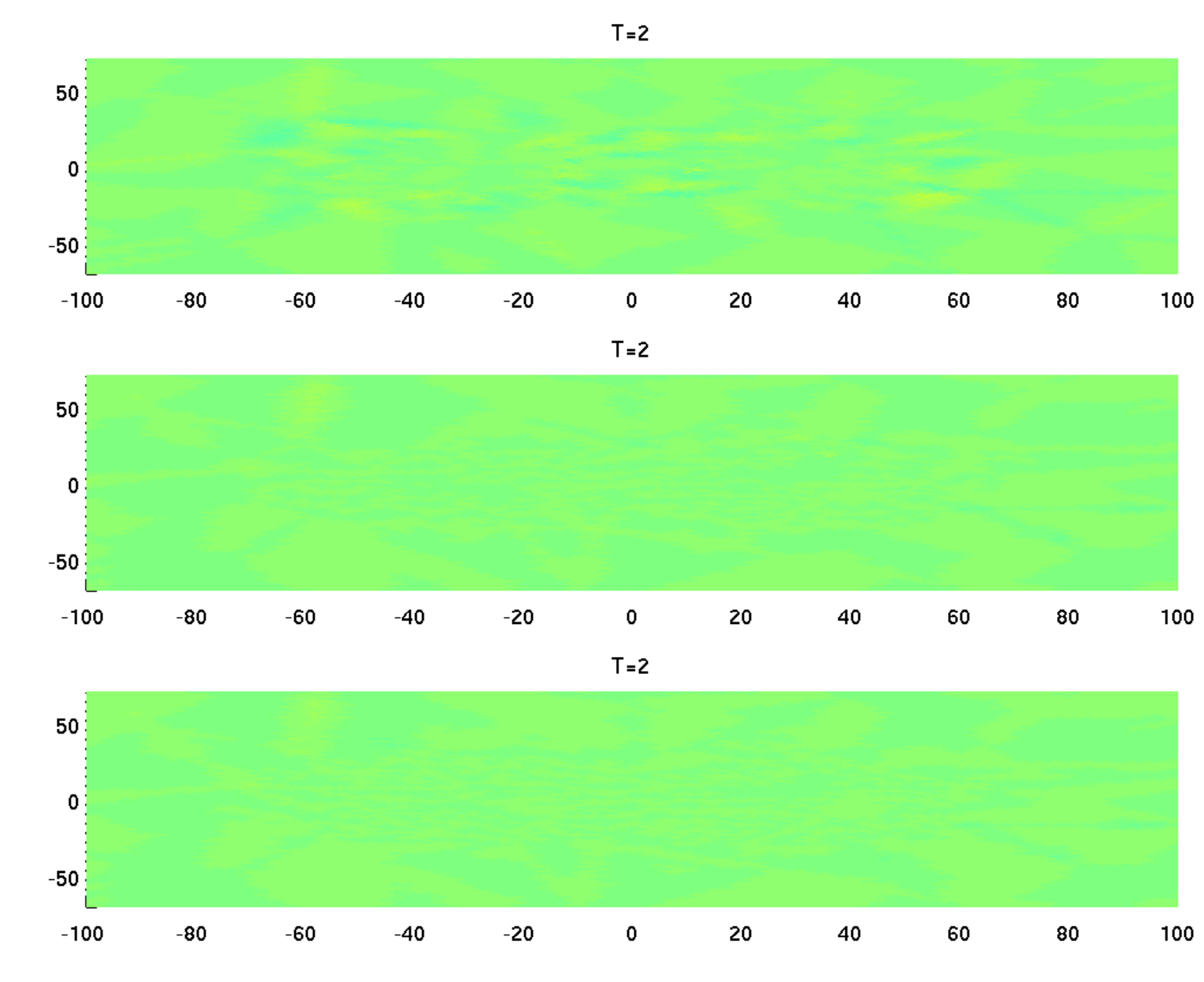}
\includegraphics[width=1.25in,height=3.5in]{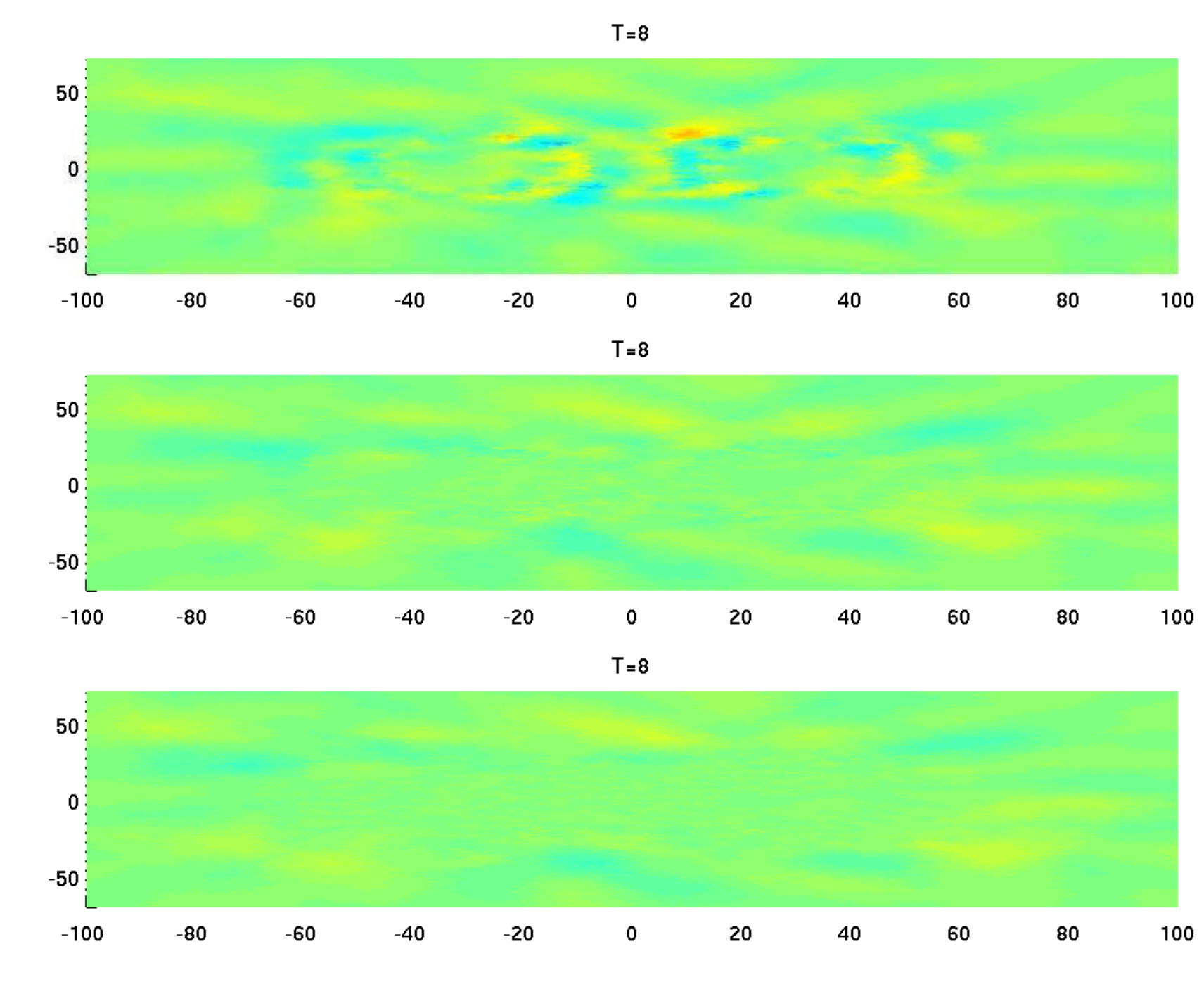}
\includegraphics[width=1.25in,height=3.5in]{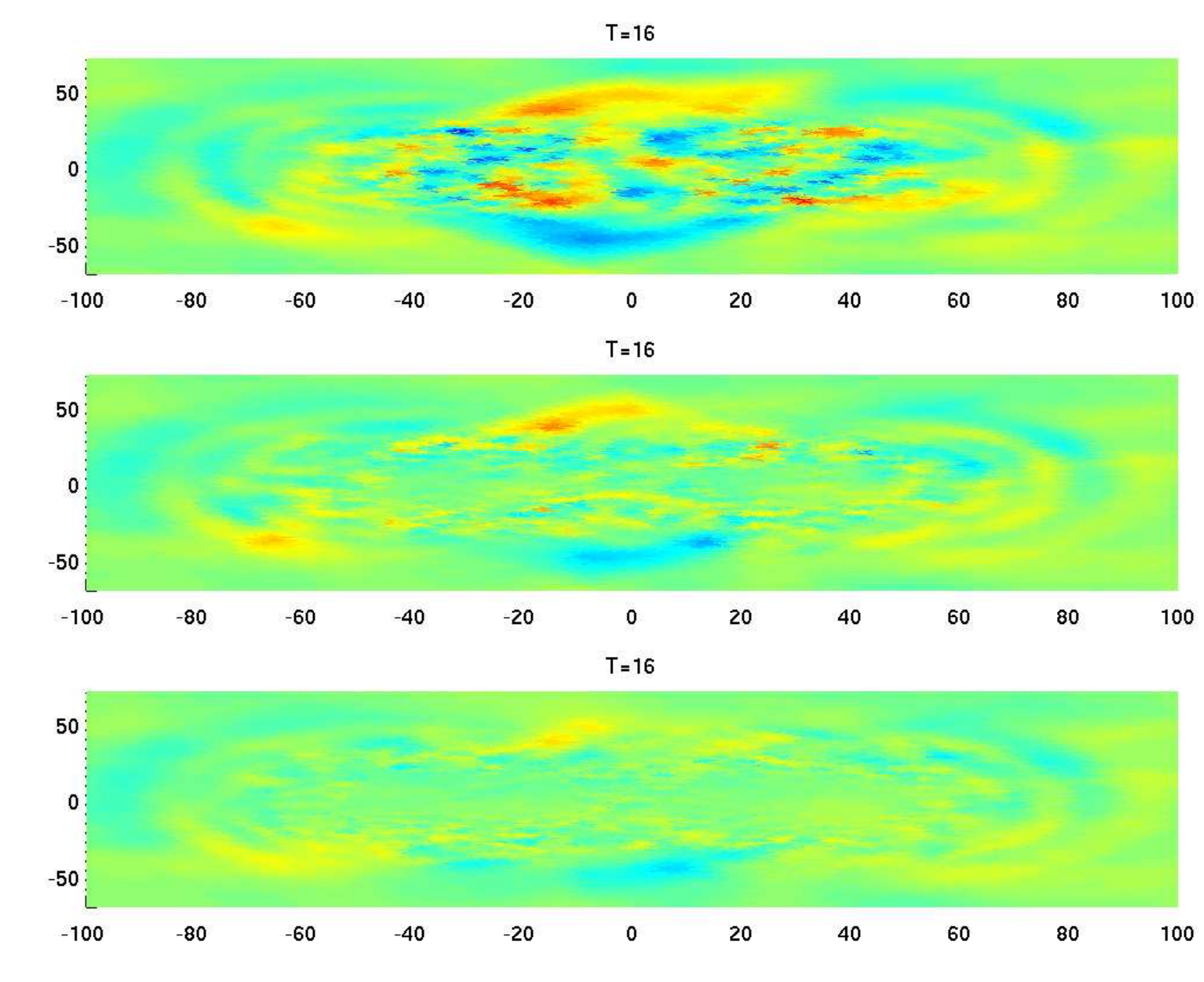}
\includegraphics[width=1.25in,height=3.5in]{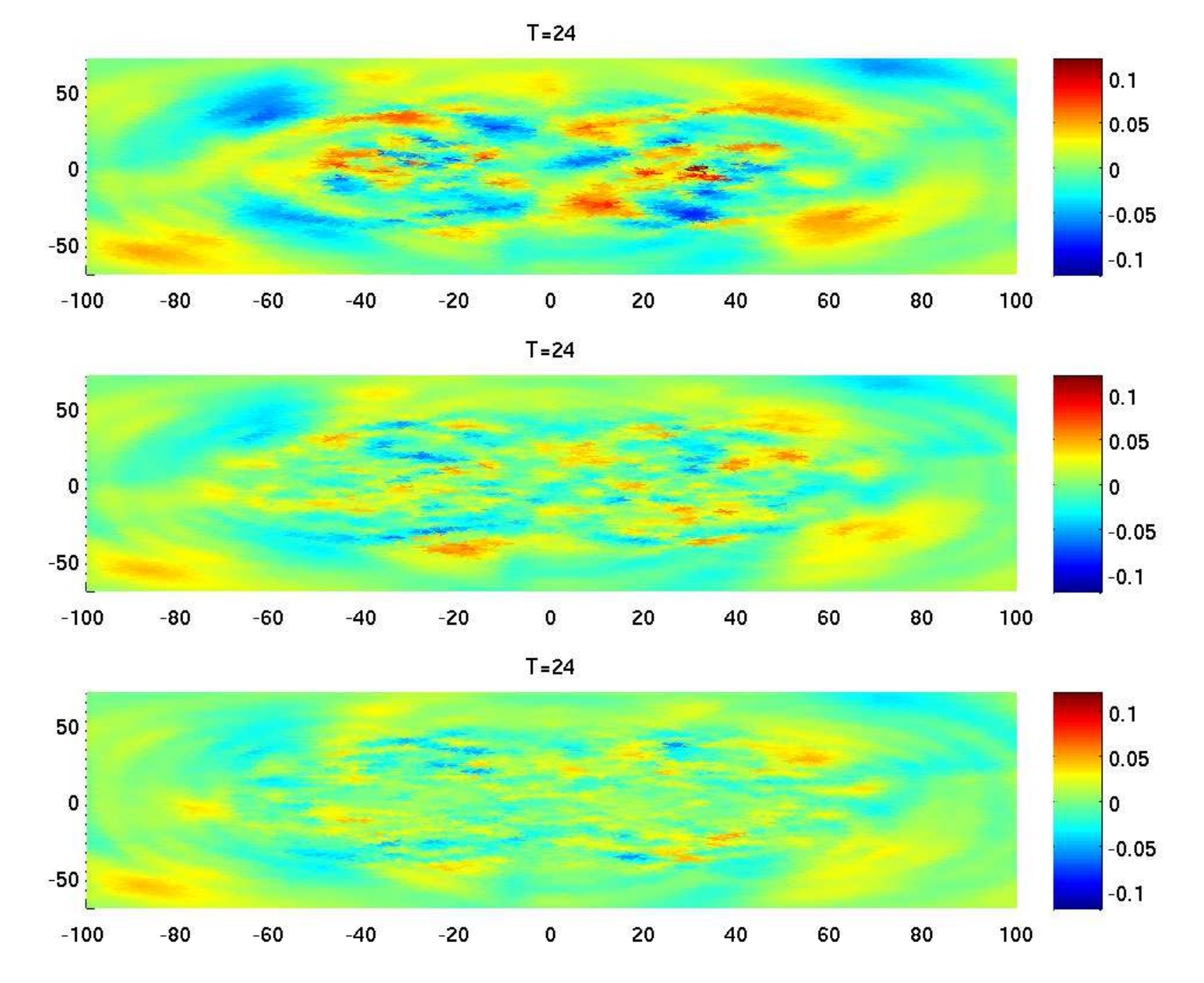}
\end{center}
\caption{Snapshots of the numerical error. From left to right, $T=2, 8, 16$, and $24$.  
From top to bottom: the error from the conventional Galerkin method, the extended Galerkin method with $L=4$, and the extended Galerkin method with $L=8$ for the Krylov subspaces \eqref{eq: krylov}. }
\label{fig: error}
\end{figure}

In Fig. \ref{fig: error}, we show the error of the Galerkin methods. Clearly, the conventional Galerkin method 
exhibits large error inside of the atomistic region, which is a result of reflection from the interface. In contrast, the extend Galerkin method introduces much less reflection, and offers better accuracy. The accuracy improves for further extensions. 

\bigskip

\section{Conclusion and discussions}
We have presented a new Galerkin formulation to coarse-grain a molecular dynamics model for crystalline solids. Our method has been designed to address the following issues:
\begin{enumerate}
\item The artificial reflections at the interface.
\item The difficulty in pre-computing the memory functions in the generalized Langevin equations.
\item The stability of CG models. 
\end{enumerate}
The method is referred to as an extended Galerkin method to distinguish it from the
conventional Galerkin method and to emphasize the variational formulation.

It is possible to extend the current formulation to finite temperature. The key is to satisfy the second
fluctuation-dissipation theorem \cite{Kubo66}. One needs to introduce a Gaussian noise to be consistent with the approximate kernel function. This is work underway.

\bibliographystyle{siam}
\bibliography{Gal0}

\end{document}